
\input  amstex
\input amsppt.sty
\magnification1200
\vsize=23.5truecm
\hsize=16.5truecm
\vcorrection{-10truemm}
\NoBlackBoxes

\def\d{d\!@!@!@!@!@!{}^{@!@!\text{\rm--}}\!}
\def\supp{\operatorname{supp}}
\def\Rn{\Bbb R^n}
\def\crp{\overline{\Bbb R}_+}

\def\crpm{\overline{\Bbb R}_\pm}

\def\rnp{{\Bbb R}^n_+}

\def\rnpm{\Bbb R^n_\pm}
\def\crnp{\overline{\Bbb R}^n_+}

\def\comega{\overline\Omega }

\def\ang#1{\langle {#1} \rangle}
\def\simto{\overset\sim\to\rightarrow}
\def\crpp{\overline {\Bbb R}^2_{++}}
\def\rpp{ {\Bbb R}^2_{++}}
\def\rp{ \Bbb R_+}
\def\rmi{ \Bbb R_-}
\define\tr{\operatorname{tr}}
\define\dist{\operatorname{dist}}

\def\wA{\widetilde A}

\def\N{\Bbb N}
\def\R{\Bbb R}
\def\C{\Bbb C}
\def\ol{\overline}
\def\weight{\ang}
\def\SD{\Cal S}
\def\F{\Cal F}
\def\Z{\Bbb Z}

\define\uu{\underline u}
\define\uv{\underline v}

\define\uB{\underline B}

\define\uomega{\underline \Omega}
\define\ucomega{\overline{\underline \Omega}}

\document

\topmatter
\title
Spectral asymptotics for nonsmooth singular Green operators
\endtitle
\author {Gerd Grubb}\endauthor
\date {With an appendix by Helmut Abels 
}\enddate
\address
{G.\ Grubb, Department of Mathematical Sciences, Copenhagen University,
Uni\-ver\-si\-tets\-par\-ken 5, DK-2100 Copenhagen, Denmark.
E-mail {\tt grubb\@math.ku.dk}}
\endaddress
\address{H.\ Abels, Faculty for
Mathematics, University Regensburg, 93040 Regensburg, Germany,
E-mail {\tt Helmut.Abels\@mathematik.uni-regensburg.de}}
\endaddress
\abstract
Singular Green operators $G$ appear typically as boundary correction
terms in resolvents for elliptic boundary value problems on a domain $\Omega \subset {\Bbb R}^n$, and more
generally they appear in the calculus of pseudodifferential boundary
problems. In particular, the
boundary term in a Krein resolvent formula is a singular Green
operator. It is well-known in
smooth cases that when $G$ is of negative order $-t$ on a bounded
domain, its eigenvalues or  $s$-numbers
have the behavior (*) $s_j(G) \sim c j^{-t/(n-1)}$ for $j\to \infty $,
governed by the boundary dimension $n-1$. In some nonsmooth cases, upper
estimates (**) $s_j(G) \le Cj^{-t/(n-1)}$ are known.

We show that (*) holds when $G$ is a general selfadjoint nonnegative
singular Green operator with symbol merely
H\"o{}lder continuous in $x$. We also show (*) with $t=2$ for the
boundary term in the Krein resolvent formula comparing the Dirichlet
and a Neumann-type problem for a strongly elliptic second-order
differential operator (not necessarily selfadjoint) with coefficients
in $W^1_q(\Omega )$ for some $q>n$. 
\endabstract
 \endkeywords
\endtopmatter

\head Introduction\endhead

Singular Green operators arise typically as boundary correction terms
in solution formulas for elliptic boundary value problems. 
For example, if $A$ is a strongly elliptic second-order differential
operator with smooth coefficients on ${\Bbb R}^n$, with inverse $Q$,
and $\Omega \subset {\Bbb R}^n$ is smooth bounded, then the solution operator for the Dirichlet problem $Au=f$ on $\Omega
$, $u|_{\partial\Omega }=0$, has the form
$$
A_\gamma ^{-1}=Q_+ +G_\gamma  ,\tag i
$$
where $Q_+$ is the truncation of the pseudodifferential operator ($\psi $do) $Q$ to $\Omega $ and $G_\gamma 
$ is a singular Green operator. Another typical singular Green
operator is the difference between the solution operators for two
different boundary value problems for $A$,
$$
G=\wA^{-1}-A_\gamma ^{-1}.\tag ii
$$
In the study of the spectral behavior it is found that whereas the
eigenvalues or singular values ($s$-numbers) of $\wA^{-1}$ and $A_\gamma ^{-1}$ have the
behavior $s_j(\wA^{-1})\sim c_Aj^{-2/n}$, the $s$-numbers of $G_\gamma
$ and $G$
behave like $s_j(G)\sim c_Gj^{-2/(n-1)}$ (cf.\ e.g.\ Grubb \cite{G74},
\cite{G84}, Birman-Solomyak \cite{BS80}). Here the boundary dimension
$n-1$ enters since these operators have their essential effect near
the boundary.

This spectral behavior of $G$ is well-known for operators $A$ with smooth
coefficients, whereas in cases with nonsmooth coefficients, generally
only upper estimates $s_j(G)\le Cj^{-2/(n-1)}$ are known (from Birman
\cite{B62}, when coefficients are continuous with bounded first
derivatives and the boundary is $C^2$).

We shall here address the question of showing asymptotic estimates for
singular Green operators in
cases with nonsmooth coefficients. The main tool will be the calculus
of nonsmooth pseudodifferential boundary operators ($\psi $dbo's)
developed by Abels \cite{A05}, \cite{A05a} (as a generalization of the
smooth $\psi $dbo's, Boutet de Monvel \cite{B71}, Grubb \cite{G84},
\cite{G96}).
Our results will deal both with general selfadjoint nonnegative
singular Green operators with $C^\tau $-smoothness in the
$x$-variable, and with the special, not necessarily selfadjoint operators
in (i), (ii), with $W^1_q$-smoothness of the coefficients of $A$ ($q>n$),
building on the resolvent construction in Abels-Grubb-Wood \cite{AGW12}.  
\smallskip

{\it Contents.} In Section 1 we recall the Krein resolvent formula in
the smooth setting and show a precise spectral asymptotic estimate in the
selfadjoint case. Section 2 deals with spectral asymptotic results for
nonsmooth $\psi $do's, recalling an early result of
Birman and Solomyak, and proving a result for $C^\tau $-smooth $\psi
$do's of negative order (any $\tau >0$), based on the calculus for such $\psi $do's by Marschall.  
In Section 3 we recall the theory of nonsmooth $\psi $dbo's by Abels, with some supplements.
Section 4 gives the proof of spectral asymptotic estimates for
nonsmooth selfadjoint singular Green operators of negative order and
class 0, defined on $\rnp$ or on a smooth bounded set $\Omega \subset{\Bbb R}^n$. In Section 5, similar results are obtained for the singular
Green term in the Dirichlet resolvent (i), and for the singular Green term
in the Krein formula (ii), by a different method based on work of
Abels, Grubb and Wood; here nonselfadjointness is allowed. In Section
6, the results of Section 5 are extended to the case of domains
$\Omega $ with $B^{\frac32}_{p,2}$-boundary; this includes
$C^{\frac32+\varepsilon }$-domains. Finally, the
Appendix written by Abels gives the proof of one of the theorems used
in Section 3.

\head {1. Results in the smooth case} \endhead

\subhead 1.1 Some notation \endsubhead

For $x\in{\Bbb R}^n$ we denote $x'=(x_1,\dots,x_{n-1})$, so that
$x=(x',x_n)$, and we denote $\rnpm=\{x\mid x_n\gtrless 0\}$. Moreover,
$\ang x=(1+|x|^2)^{\frac12}$.

When $\Omega $ is a smooth open subset of ${\Bbb R}^n$ with boundary $\partial\Omega =\Sigma $, we use the standard
$L_2$-Sobolev spaces, with the following notation: $H^s({\Bbb R}^n)$
($s\in 
{\Bbb R}$) has the norm
$\|v\|_s
=\|{\Cal F}^{-1}(\ang\xi ^s{\Cal
F}v)\|_{L_2({\Bbb R}^n)}$; here ${\Cal F}$ is the Fourier transform and
$\ang\xi =(1+|\xi |^2)^{\frac12}$. Next, $H^s(\Omega )=r_{\Omega }H^s({\Bbb
R}^n)$ where $r_\Omega $ restricts to $\Omega $, provided with
the norm $\|u\|_{s}=\inf\{\|v\|_s\mid v\in H^s({\Bbb
R}^n),\, u=r_\Omega v \}$. Moreover, we denote by   $H^s_0(\comega )$
the space $\{u\in H^s({\Bbb R}^n)\mid \operatorname{supp}u\subset \comega \}$;
it is a closed subspace of $H^s({\Bbb R}^n)$, and there is an
identification of the antidual space of $H^s(\Omega )$ (the space
of antilinear, i.e., conjugate linear, functionals), with
$H^{-s}_0(\comega)$ for any $s\in{\Bbb R}$, with a duality consistent
with the $L_2(\Omega )$ scalar product. Spaces over the boundary,
$H^s(\Sigma )$, are defined by local coordinates from $H^s({\Bbb
R}^{n-1})$, $s\in{\Bbb R}$. 
Here there is an identification of
$H^{-s}(\Sigma )$ with the antidual space of $H^s(\Sigma )$.

Occasionally we shall also refer to some $L_q$-based Sobolev spaces, $1<q<\infty $,
namely the Bessel-potential spaces $H^s_q(\Omega )$; for $s=k\in{\Bbb
N}_0$ they are also denoted $W^k_q(\Omega )$. They are defined in a
similar way as
above from
$H^s_q({\Bbb R}^n)$, provided with the norm $\|v\|_{s,p}
=\|{\Cal F}^{-1}(\ang\xi ^s{\Cal
F}v)\|_{L_q({\Bbb R}^n)}$. One has for $k\in{\Bbb N}_0$ that
$H^k_q(\Omega )=W^k_q(\Omega )=\{u\in L_q(\Omega )\mid D^\alpha u\in L_q(\Omega
)\text{ for }|\alpha |\le k\}$ (allowing also $q=\infty $). 

For $\tau \ge 0$,  we denote by $C^r(\Omega )$ the space of 
continuous 
functions $f(x)$ such
 that when $[\tau ]$ is the largest integer $\le \tau $,
$$
\|f\|_{C^\tau }\equiv \sum_{|\beta |\le [\tau ]}\sup _x|D^\beta f(x)| +
\sum_{|\beta |= [\tau ]}\sup _{x\ne y}\frac{|D^\beta f(x)-D^\beta f(y)|}{|x-y|^{\tau -[\tau ]}}<\infty .\tag1.1
$$
For  $\tau =k\in{\Bbb N}_0$, they are also denoted $C^k_b(\Omega )$; when
$\tau =k+\sigma $ with $k\in{\Bbb N}_0$ and $\sigma \in\,]0,1[\,$,
they are the H\"o{}lder spaces $C^{k,\sigma }(\Omega )$.
  We denote $\bigcap_{\tau \ge 0}C^\tau = C^\infty _b$.

\subhead 1.2 The Krein resolvent formula \endsubhead

As a point of departure, consider the Krein resolvent formula
comparing the resolvents of the Dirichlet realization $A_\gamma $ and
a Neumann-type realization $\wA$. The operators are defined from a second-order strongly elliptic
operator $A$ with smooth coefficients on a smooth bounded or exterior
domain $\Omega $ with boundary $\Sigma $ (cf.\ Section 5) by the boundary conditions,
assumed elliptic,
$$
\gamma _0u=0,\text{ resp. }\chi u=C\gamma _0u,\text{ on }\Sigma ,
$$
with a first-order tangential differential operator $C$; here $\gamma
_0u=u|_{\Sigma }$ and  $\chi u$ is the conormal derivative (5.12). Then one has when $0\in \varrho (A_\gamma )\cap \varrho (\wA)$ (the
resolvent sets) that
$$
\wA^{-1}-A_\gamma ^{-1}= K_\gamma L^{-1}(K'_\gamma)^*,\tag1.2
$$
where $K_\gamma $ and $K'_\gamma $ are the Poisson operators for the
Dirichlet problem for $A$ resp.\ $A'$, $L$ is the realization of
the first-order elliptic $\psi $do $C-P_{\gamma ,\chi }$ with
$D(L)=H^{\frac32}(\Sigma )$, and $P_{\gamma ,\chi} =\chi K_\gamma $,
the Dirichlet-to-Neumann operator. Cf.\ Grubb \cite{G68}, \cite{G74},
and formulas (2.10), (2.15), (3.45) in  Brown-Grubb-Wood \cite{BGW09}
(based on \cite{G68}, \cite{G74}); see also e.g.\ \cite{G08}. Abstract
versions have been known for many years and various concrete applications to
elliptic PDE given recently, cf.\ e.g.\  Malamud
\cite{M10}, Gesztesy-Mitrea \cite{GM11}, and their listings of other contributions.

Consider for precision in this introductory section the case where the operators $A_\gamma $ and $\wA$ are
selfadjoint; then so is $L$, as an operator from $H^{-\frac12}(\Sigma
)$ to $H^\frac12(\Sigma )$. Then also $K_\gamma =K'_\gamma
$. (Nonselfadjoint cases, where principal estimates of $s$-numbers can
be obtained, are treated in \cite{G12}, Sect.\ 10.) 

The operator $G=K_\gamma L^{-1}K_\gamma ^*$ is a singular Green
operator of order $-2$, and its spectral asymptotic behavior can be found in the
following way: Denote the positive resp.\ negative eigenvalues by $\mu
^\pm_j$, monotonely ordered and repeated according to multiplicity
(one of the sequences may be finite and then needs no treatment, we
leave this aspect out of the explanation). Then we have, since $\mu _j^\pm(B_1B_2)=\mu _j^\pm(B_2B_1)$:
$$
\mu _j^\pm(G)=\mu _j^\pm(L^{-1}K_\gamma ^*K_\gamma )=\mu _j^\pm(L^{-1}P_1),\tag1.3
$$
 where $P_1=K^*_\gamma K_\gamma $ is a selfadjoint positive $\psi $do
 on $\Sigma $ of order $-1$ (it is also used in \cite{G11} Th.\
 3.4, which gives more information).  Let $P_2=P_1^{\frac12}$, positive selfadjoint of order $-\frac12$, then
 we have furthermore:
$$
\mu _j^\pm(G)=\mu _j^\pm(L^{-1}P_2^2)=\mu _j^\pm(P_2L^{-1}P_2)=\mu
_j^\pm(P_3),\quad P_3=P_2L^{-1}P_2.\tag1.4
$$
Estimates of the eigenvalues $\mu _j^\pm$ are connected in a known way
with estimates of the corresponding counting functions
$N^{\prime\pm}(t;G)=N^{\prime\pm}(t;P_3)$;  here $N^{\prime\pm}(t;S)$ 
indicates the number of positive, resp.\ negative
eigenvalues of $S$ outside the interval $\,]-1/t,1/t[\,$. 

For the counting functions we apply the results of H\"o{}rmander \cite{H68} and Ivrii \cite{I82} to the selfadjoint elliptic
$\psi $do $P_3^{-1}$ of order 2. This gives:

\proclaim{Theorem 1.1} In the smooth selfadjoint case one has for the
operator $\wA^{-1}-A_\gamma ^{-1}=K_\gamma L^{-1}K_\gamma ^*=G$ that the eigenvalues satisfy
$$ 
N^{\prime\pm}(t;\wA^{-1}-A_\gamma ^{-1})=C ^\pm
t^{(n-1)/2}+O(t^{(n-2)/2})\text{ for }t\to\infty ,\tag1.5
$$
where $C^{\pm}$ are determined from the principal symbol $p_3^0(x',\xi ')$ of $P_3$
(defined through {\rm (1.3)--(1.4)}).
Moreover, if $(p_3^0)^{-1}$ 
satifies  Ivrii's
  condition from {\rm \cite{I82}} (the
  bicharacteristics through points of $T^*(\Sigma )\setminus 0$ are
  nonperiodic except for a set of measure zero), there are
  constants $C_1^\pm$ such that 
$$
N^{\prime\pm}(t;\wA^{-1}-A_\gamma ^{-1})=C^\pm t^{(n-1)/2}+C_1 ^\pm t^{(n-2)/2}
+o(t^{(n-2)/2})\text{ for }t\to\infty .\tag1.6
$$
\endproclaim

The estimate (1.5) follows from H\"o{}rmander \cite{H68} in the
scalar case, Ivrii allows abitrary elliptic systems and has the
precision in (1.6). (The result does not seem to have been formulated
with this precision before.)

Note that in the scalar case, if $\Sigma $ is connected, the
selfadjointness and ellipticity prevents the
principal symbol of $L$ from changing between positive and negative
values, since it must be real. So in that case only one of the
sequences $\mu _j^+$ or $\mu _j^-$ is infinite.

In rough cases, we do not expect to get the fine remainder estimates, 
but will aim for principal asymptotic (Weyl-type) estimates, as
obtained in \cite{G84} for general singular Green operators in the
smooth case.

\head 2. Spectral estimates for nonsmooth pseudodifferential operators  
\endhead

\subhead 2.1  Weak Schatten classes \endsubhead

As in \cite{G84} we denote by ${\frak C}_p(H,H_1)$ the  $p$-th Schatten class 
consisting of the compact operators $B$ from a Hilbert space $H$ to
another $H_1$ such that
$(s_j(B))_{j\in{\Bbb N}}\in \ell_p({\Bbb N})$. Here $s_j(B)=\mu
_j(B^*B)^\frac12$, where $\mu _j(B^*B)$ denotes the $j$-th positive
eigenvalue of $B^*B$, arranged nonincreasingly and repeated according
to multiplicities. 
The so-called weak Schatten class
consists of the compact operators $B$ such that
$$
s_j(B)\le Cj^{-1/p}\text{ for all }j;\text{ we set } {\bold N}_p(B)=\sup_{j\in{\Bbb N}}s_j(B)j^{1/p};\tag2.1
$$
the notation $\frak S_{(p)}(H,H_1)$ was used in
\cite{G84} for this space.
(The indication $(H,H_1)$ is replaced by $(H)$ if $H=H_1$; it can be omitted when it is clear from the context.)
Different notation is used in some other
works; $\frak C_p$ is sometimes called $\frak S_p$, and
$\frak S_{(p)}$ is also sometimes called $\Sigma _p$, $\frak S_p$ or $\frak
S_{p,\infty }$. To avoid confusion, we shall in the present paper use
the notation $\frak S_{p,\infty }$ for the $p$-th weak Schatten class. 

We recall (cf.\ e.g.\ \cite{G84} for details and references) that ${\bold N}_p(B)$ is
a quasinorm on $\frak S_{p,\infty }$, satisfying
$$
\aligned
{\bold N}_p(\sum_{k=1}^{k_0}B_k)&\le C_{p,\gamma }\sum_{k=1}^{k_0}{\bold N}_p(B_k)
k^\gamma ,\\\text{ with }
\gamma =0\text{ if }&p>1,\quad \gamma >p^{-1}-1\text{ if }p\le 1;
\endaligned\tag2.2$$
here $C_{p,\gamma }$ is independent of $k_0$. Recall also that
$$
\frak S_{p,\infty }\cdot \frak S_{q,\infty }\subset \frak S_{r,\infty }, \text{ where }r^{-1}=p^{-1}+q^{-1},\tag2.3
$$
and
$$
s_j(B^*)=s_j(B),\quad
s_j(EBF)\le \|E\|s_j(B)\|F\|,\tag2.4
$$
when $E\colon H_2\to H$ and $F\colon H_1\to H_3$ are bounded linear maps between Hilbert
spaces. 

Moreover, we recall that when $\Xi $ is a bounded open subset
of ${\Bbb R}^m$ and reasonably regular, then the injection $H^t(\Xi
)\hookrightarrow L_2(\Xi )$ is in $\frak S_{m/t,\infty }$ when $t>0$. It
follows that when $B$ is a linear operator in $L_2(\Xi )$ that is
bounded from $L_2(\Xi )$ to $H^t(\Xi )$, then $B\in \frak
S_{m/t,\infty }$, with
$$
{\bold N}_{m/t}(B)\le C\|B\|_{\Cal L(L_2(\Xi ),H^t(\Xi ))}.\tag2.5
$$

\subhead 2.2  Results of Birman and Solomyak \endsubhead
 
We shall study pseudodifferential operators 
with $C^\tau $-smoothness in the
$x$-variable further below, but let us first consider some results of
Birman and Solomyak \cite{BS77}.
They show an asymptotic  result under
weak 
smoothness hypotheses both in the $x$- and the $\xi $-variable:

\proclaim{Theorem 2.1}{\rm \cite{BS77}} On a closed manifold $\Xi $ of dimension
$m$, let $P$ be defined in local
coordinates from symbols $p(x,\xi ) $ that are homogeneous in $\xi $
of degree $-t<0$. Denote $m/t=\mu $.

$1^\circ$ One has that $P\in \frak S_{\mu ,\infty }$ (i.e., $s_j(P)$ is
$O(j^{-t/m})$) under the following
hypotheses on the symbol in local coordinates: 

If $\mu \le 1$ (i.e.,  $t\ge m$), assume that the symbols restricted to $\xi \in S^{m-1}=\{|\xi |=1\}$ are in 
$ L_\infty (S_\xi ^{m-1}, C_x^\varepsilon)$ for some $\varepsilon >0$.

If $\mu >1$ (i.e., $t<m$) assume that the symbols at $\xi \in S^{m-1}$ are in 
$ L_\infty (S ^{m-1}, W^l_p)$, where 
$$
p\ge 2, \quad pl>m,\quad 1/p>1/2-1/q_1,\tag2.6 
$$
for some $2<q_1,q_2\le \infty $ with $q_1^{-1}+q_2^{-1}=\mu ^{-1}$, $l\in\rp$.

$2^\circ$ There is an asymptotic estimate 
$$
s_j(P)j^{t/m}\to c,\tag2.7
$$
when the properties under $1^\circ$ hold with $L_\infty $ replaced by $C^0$. 
\endproclaim

The constant $c$ is the same as described further below in Theorem 2.5.

It is particularly interesting here that these estimates allow nonsmoothness in
$\xi $, namely just boundedness or continuity, with no requirement on $\xi
$-derivatives.

Also some nonsmoothness in $x$ is allowed, best in low dimensions. Let
us see what it means for $t=2$. 

Here $1^\circ$ applies when
when $\mu \le 1$, i.e., $m\le 2$, so only $C^\varepsilon $-smoothness
in $x$ is needed then.

When $2<m\le 4$, so that $1<\mu \le 2$, $1/2\le 1/\mu <1$, we can
in $2^\circ$ take $1/q_1=1/2-\delta $ with a small positive $\delta $ (hereby $q_1>2$) and
$1/q_2=1/\mu -1/q_1=1/\mu -1/2+\delta  $  (which is $<1/2$ for
$\delta  <1-1/\mu $, hereby $q_2>2$). Then the requirement
$1/p>1/2-1/q_1 $ allows taking $1/p$ arbitrarily close to 0.
Now $pl>m$ can be fulfilled for arbitrarily small $l$ by taking $p$
sufficiently large, so $W^l_p$ can be taken to contain $C^\varepsilon $,
for a given small $\varepsilon $. 

When $m>4$, the inequalities will put a positive lower limit on the
possible $l$ that can enter in (2.6). For example, for $m=5$, $\mu
=5/2$ and $1/\mu =2/5$, then we can at best take $q_1=5/2$ and $q_2=\infty $,
which restricts $p$ by 
$$
1/p>1/2-1/q=1/2-2/5=1/10,
$$
so that $p<10$. Then $pl>m$ is at best obtained with $l>1/2$.

So already for $m=5$, hence for the boundary of a set of dimension $6$, the Birman-Solomyak result will not give
asymptotic estimates for the most general situation in
\cite{AGW12} where the symbols are only $C^\tau $ in $x$ with a $\tau
<1/2$.

It should be noted that our symbol classes have a high degree of smoothness in $\xi $ in
contrast to those of Birman and Solomyak; we do not need their generality for the
present purposes, and have
to find another point of view.

\subhead {2.3 Approximation of $\psi $do's by operators with smooth symbols}\endsubhead

We now turn to the symbol classes with H\"o{}lder smoothness in $x$
and full smoothness in $\xi $, as defined in Kumano-go and Nagase
\cite{KN78}, Marschall \cite{M87}, Taylor \cite{T91} and other places:

\proclaim{Definition 2.2} Let $d\in{\Bbb R}$, $\tau >0$, $m\in{\Bbb
N}$, $N\in {\Bbb N_0}$.
The space $C^\tau S^{d}_{1,0}({\Bbb R}^m\times{\Bbb R}^m,N)$ consists
of the
functions $p(x,\xi )$  of $x,\xi \in{\Bbb R}^m$ such that
$$
\|D_\xi ^\alpha p(x,\xi )\|_{C^\tau }\text{ is }O(\ang{\xi}^{d-|\alpha
|})\text{ for }|\alpha |\le N.\tag2.8 
$$
We denote  $\bigcap_{N\in{\Bbb N}_0}C^\tau S^{d}_{1,0}({\Bbb
R}^m\times{\Bbb R}^m,N)=C^\tau S^{d}_{1,0}({\Bbb R}^m\times{\Bbb
R}^m)$.

The symbol $p(x,\xi )$ is said to be {\bf polyhomogeneous} (with step $1$), when there is an
asymptotic expansion in symbols $p_{d-j}(x,\xi )$ homogeneous of order
$d-j$ in $\xi $ for $|\xi |\ge 1$, in the sense that each  $p_{d-j}\in C^\tau
S^{d-j}_{1,0}({\Bbb R}^m\times{\Bbb R}^m,N)$, and for all $J$,
$p-\sum_{j<J}p_{d-j}$ is in $C^\tau S^{d-J}_{1,0}({\Bbb
R}^m\times{\Bbb R}^m,N)$. For the subspaces of polyhomogeneous symbols
we use the notation $C^\tau S^{d}$ instead of $C^\tau S^{d}_{1,0}$.
\endproclaim

Our convention for the Fourier transform is: ${\Cal F}_{x\to\xi }
u=\int_{{\Bbb R}^m}e^{-ix\cdot\xi }u(x)\,dx$; then  ${\Cal
F}^{-1}_{\xi \to x }f=\int_{{\Bbb R}^m}e^{ix\cdot\xi }f(\xi )\,\d\xi$,
with $\d\xi =(2\pi )^{-n}d\xi $.

A symbol in $C^\tau S^d_{1,0}({\Bbb
R}^m\times{\Bbb R}^m)$ defines a $\psi $do $P$ by 
$$
Pu=\operatorname{OP}(p(x,\xi ))u=\int e^{i(x-y)\cdot
\xi }p(x,\xi ) u(y)\,dy\d\xi , \text{ also called }p(x,D_x)u,\tag2.9
$$
with a suitable interpretation of the integral (as an oscillatory integral). 
We recall
in passing that one can also define operators ``in $y$-form'' resp.\ ``in
$(x,y)$-form'' by formulas
$$
\aligned
\operatorname{OP}(p(y,\xi ))u&=\int e^{i(x-y)\cdot
\xi }p(y,\xi ) u(y)\,dy\d\xi,\text{ resp.\ }\\
 \operatorname{OP}(p(x,y,\xi ))u&=\int e^{i(x-y)\cdot
\xi }p(x,y,\xi ) u(y)\,dy\d\xi ;
\endaligned\tag2.10
$$
the latter is in some texts said to be defined from a double
symbol. We say that (2.9) is ``in $x$-form''.

It is well-known that $P$ in (2.9) satisfies
$$
P\colon H^{d}({\Bbb R}^m )\to L_2({\Bbb R}^m ) \text{ and hence }P^*\colon
L_2({\Bbb R}^m )\to H^{-d}({\Bbb R}^m ).\tag2.11
$$
(More information is given below in Theorem 2.3.)
Hence when $d$ is a negative number $-t$, and $P$ is defined on a
compact $m$-dimensional manifold from such symbols in local
coordinates, $P\in \frak S_{m/t,\infty }$ in view of (2.5). 

This gives upper spectral estimates, and we shall obtain the asymptotic spectral estimates by approximation of the $\psi $do
symbols by smooth polyhomogeneous symbols for which the estimates are
known. The so-called {\it symbol smoothing}, where $p$ is written as
$p^\sharp+p^{\flat}$, $p^\sharp\in S^m_{1,\delta }$ and $p^{\flat}$ of lower order,
is not useful here, since the polyhomogeneity is lost in this
decomposition. What we do is in fact simpler, namely approximation
by convolution in the $x$-variable with
an approximate unit.

There is one small obstacle here, namely that $C^\infty $-functions
are not dense in $C^\tau $ when $\tau $ is not integer. 
For example, if one compares a $C^\infty $-function $u (t)$ on
$[-1,1]$ with the function $|t|^\tau \in C^\tau ([-1,1])$ for some
$\tau \in\,]0,1[\,$, one finds that 
$$
\sup _{t\ne s}\frac{|(u (t)-|t|^\tau)
-(u (s)-|s|^\tau )|}{ |t-s|^{\tau }}
\ge \sup _{t\ne 0}\frac{|u (t)
-u (0)-|t|^\tau|}{ |t|^{\tau }}
\ge 1.
$$
However, it is well-known (cf.\ e.g.\ Lunardi \cite{L95} Ch.\ 1) that the so-called little-H\"ol{}der space $h^\tau ({\Bbb R}^m)$,
consisting of functions $u(x)\in C^\tau ({\Bbb R}^m)$ such that  
$$
\lim_{h\to 0}\sup_{0<|x-y|\le h}\frac{|u(x)-u(y)|}{|x-y|^\tau }=0,\tag2.12
$$
is a closed subspace of $C^\tau ({\Bbb R}^m)$ for $0<\tau
<1$ that equals the closure of $C^\infty _b({\Bbb R}^m)$  in the
$C^\tau $-norm. Here, when $\varrho _k $ is an approximate unit, i.e.,
$\varrho  _k(x)=k^m\varrho (kx )$ for $k\in{\Bbb N}$,
for some $\varrho \in C_0^\infty ({\Bbb R}^m)$ with $\|\varrho
\|_{L_1({\Bbb R}^m)}=1$, one can check that if 
$u\in h^\tau ({\Bbb R}^m)$, then $\varrho _k*u\to u$ in $h^\tau
({\Bbb R}^m)$ for $k\to \infty $.

Now for $0<\tau <\tau _1<1$,
$$
C^{\tau _1}({\Bbb R}^m)\hookrightarrow h^\tau ({\Bbb
R}^m)\hookrightarrow C^\tau ({\Bbb R}^m),\tag2.13
$$
so a function $u\in C^{\tau _1}({\Bbb R}^m)$ is approximated in $C^\tau
$-norm by $\varrho _k*u$ for $k\to \infty $.

Similarly, one can check that when 
$p(x,\xi )$ is a symbol in $ C^{\tau } S^d_{1,0}({\Bbb R}^m\times{\Bbb
R}^m,N)$, then
$$
\varrho _k(x) * p(x,\xi )\to p(x,\xi ) \text{ in } C^{\tau '}S^d_{1,0}({\Bbb
R}^m\times{\Bbb R}^m,N),\text{ when }0<\tau '<\tau <1;\tag2.14
$$
here $\varrho  
_k*p\in C^\infty _b 
S^d_{1,0}({\Bbb R}^m\times{\Bbb R}^m,N)$.

There are analogous definitions of symbol spaces where $C^\tau $ is
replaced by $H^r_q$-spaces (Bessel-potential spaces), that we shall
refer to in Section 5 below; $\psi $do's defined from such symbols
were studied by Marschall in \cite{M88}. 
 When 
$p(x,\xi )\in H^r_q S^d_{1,0}({\Bbb R}^m\times{\Bbb
R}^m,N)$, then $\varrho  
_k*p\in C^\infty _b 
S^d_{1,0}({\Bbb R}^m\times{\Bbb R}^m,N)$, and
$$
\varrho _k(x) * p(x,\xi )\to p(x,\xi ) \text{ in } H^r_qS^d_{1,0}({\Bbb
R}^m\times{\Bbb R}^m,N).\tag2.15
$$

One can also define symbols valued in Banach spaces. Let $X$ be a
Banach space, then $C^\tau S^{d}_{1,0}({\Bbb R}^m\times{\Bbb R}^m,N;X)$ consists
of the
functions $p(x,\xi )$  from  $(x,\xi )\in{\Bbb R}^{m}\times{\Bbb R}^m$ to $X$ such that
$$
\|D_\xi ^\alpha p(x,\xi )\|_{C^\tau ({\Bbb R}^{m}_x,X)}\text{ is }O(\ang{\xi}^{d-|\alpha
|})\text{ for }|\alpha |\le N,\tag2.16
$$ 
where $C^\tau ({\Bbb R}^{m}_x,X)$ is provided with the norm in (1.1)
with absolute values replaced by $X$-norms. $X$ can in particular be a
space of bounded linear operators $X=\Cal L(X_0,X_1)$ between Banach
spaces $X_0,X_1$. The use of special Fr\'e{}chet spaces such as  $\Cal
S_{+}$ in the place of $X$ is discussed below in Section 3.

The following result was shown (in a greater generality) in Marschall \cite{M87} Th.\ 2.1:

\proclaim{Theorem 2.3}
For $d\in{\Bbb R}$, $|s|<\tau $, one has that when
$p(x,\xi )\in C^\tau S^{d}_{1,0}({\Bbb R}^m\times {\Bbb R}^m,N)$ with
$N>m/2+1$, then $\operatorname{OP}(p)$ is continuous:
$$
\|\operatorname{OP}(p)\|_{\Cal L(H^{s+d}({\Bbb R}^m),H^s({\Bbb R}^m))}< \infty .
\tag2.17
$$

\endproclaim

In other words, the linear map OP  from the Banach space $C^\tau
S^{d}_{1,0}({\Bbb R}^m\times {\Bbb R}^m,N)$ to the Banach space 
$\Cal L(H^{s+d}({\Bbb
R}^m),H^s({\Bbb R}^m))$ is bounded for each $|s|<\tau $, when $N>m/2+1$.
Then also 
$$
\|\operatorname{OP}(p)-\operatorname{OP}(p_k)\|_{\Cal L(H^{s+d}({\Bbb
R}^m),H^s({\Bbb R}^m))}\to 0 \text{ for }k\to \infty ,\text{ when
} p_k=\varrho _k*p .\tag2.18
$$
This holds for each $|s|<\tau $, since there is room for a $\tau '\in
\,]0,\tau [\,$ such that $|s|<\tau '$, and the symbol convergence
holds in $C^{\tau '}S^d_{1,0}$.

To show spectral 
asymptotic estimates, we shall use what is known in smooth cases and
extend it to nonsmooth cases by use of suitable perturbation results for
$s$-numbers:

\proclaim{Lemma 2.4}

$1^\circ$ If $s_j(B)j^{1/p}\to C_0$ and $s_j(B')j^{1/p}\to 0$ for
$j\to\infty$, then $s_j(B+B')j^{1/p}\to C_0$ for $j\to\infty$.

$2^\circ$ If $B=B_M+B'_M$ for each $M\in\Bbb N$, where
$s_j(B_M)j^{1/p}\to C_M$ for $j\to\infty$ and $s_j(B'_M)j^{1/p}\le
c_M$ for
$j\in\Bbb N$, with $C_M\to C_0$ and $c_M\to 0$ for
$M\to\infty$, then $s_j(B)j^{1/p}\to C_0$ for $j\to\infty$.
\endproclaim

The statement in $1^\circ$ is the Weyl-Ky Fan theorem (cf.\ e.g.\
\cite{GK69} Th.\ II 2.3), and $2^\circ$ is a refinement shown in
\cite{G84} Lemma 4.2.$2^\circ$.

We have the following result for nonsmooth $\psi $do's on closed
smooth manifolds.

\proclaim{Theorem 2.5} Let $E$ be an $M$-dimensional smooth vector bundle $E$ over a smooth compact
boundaryless $m$-dimensional manifold $\Xi $. Let $t>0$ and 
 $0<\tau <1$, and let $P$ be a $C^\tau $-smooth
$\psi $do acting in $E$,  
with symbol defined in local
trivializations from symbols $p(x,\xi ) $ in 
$C^\tau
S^{-t}({\Bbb R}^m\times{\Bbb R}^m)\otimes \Cal L({\Bbb C}^M)$. Then the $s$-numbers of $P$ satisfy
the asymptotic estimate 
$$
s_j(P)j^{t/m}\to c(p^0)^{t/m} \text{ for }j\to\infty ,\tag2.19
$$
where 
$$
 c(p^0)=\tfrac1{m (2\pi )^m}\int_{\Xi}\int_{|\xi |=1}\tr \bigl((p^0(x,\xi )^*p^0(x,\xi ))^{m/2t}\bigr)\, d\omega dx.\tag2.20
$$
\endproclaim

\demo{Proof} 
Using the approximation (2.14) in localizations, we can
approximate $P$ by a sequence of operators $P_k$ with polyhomogeneous $C^\infty $-symbols 
(locally in  $C^\infty 
S^{-t}$), converging in the topology of symbols
in $C^{\tau '}
S^{-t}$ for $0<\tau '<\tau $.
Then the norm of $P-P_k$ in $\Cal L(H^{-t}(\Xi ), L_2(\Xi ))$
goes to $0$ for $k\to\infty $. The statement of the theorem holds for
the $P_k$, as a corollary to Seeley \cite{S67}, cf.\ \cite{G84}, Lemma
4.5ff. Moreover, since $(P-P_k)^*\to 0$ in $\Cal L(L_2(\Xi ),H^{t}(\Xi ))$,
$$
\sup_js_j(P-P_k)j^{t/m}\to 0 \text{ for }k\to\infty ,
$$
cf.\ (2.5). We also have that $c(p^0_k)\to c(p^0) \text{ for }k\to\infty $, since
the symbol sequence converges in $C^{\tau '}$. Then the conclusion follows for $P$ by Lemma 2.4 $2^\circ$.\qed

\enddemo

\head 3. Nonsmooth pseudodifferential boundary operators  \endhead

\subhead 3.1 Boundary symbols with H\"o{}lder smoothness \endsubhead 

We want to generalize the results of \cite{G84} on singular Green
operators to nonsmooth symbols. This takes a larger effort, since the
continuity from $H^{-t}(\Omega )$ to $L_2(\Omega )$, when $G$ is of order
$-t$ and class 0 on the manifold $\comega $ with boundary, only implies
$G\in \frak S_{n/t,\infty }$, not $G\in \frak S_{(n-1)/t,\infty }$ as is known in
the smooth case. 

We denote by  $\Cal S_\pm=\Cal S(\crpm)$ the space of restrictions
$r^\pm$ to $\crpm$ of functions in $\Cal S({\Bbb R})$, and by $\Cal S_{++}=\Cal S(\crpp)$ 
the space of restrictions to
$\crpp=\crp\times\crp$ of functions in $\Cal S({\Bbb R}^2)$; 
here $\Cal S({\Bbb R}^n)$ is the Schwartz space of rapidly decreasing
$C^\infty $-functions.

Singular Green operators of class 0 are defined in the smooth case on
$\crnp$ from functions $\tilde g(x',x_n,y_n,\xi ')\in S^d_{1,0}({\Bbb
R}^{n-1}\times {\Bbb R}^{n-1},\Cal S_{++})$ by
$$
Gu=\int_{{\Bbb R}^{n-1}}e^{ix'\cdot\xi '}\int_0^\infty \tilde g(x',x_n,y_n,\xi '){\Cal F}_{y'\to\xi '}u(y',y_n)\,dy_n\d\xi '.
$$
The tool that we shall use to invoke the effect of the boundary
dimension is a decomposition of the spaces $\Cal S_+$ and $\Cal
S_{++}$ in terms of (a variant of) Laguerre functions, using expansions in
such functions not only for the symbols but also on the operator
level (Section 3.3).

Let us briefly recall the Laguerre expansions (introduced in
\cite{B71}, their role being further explained in \cite{G84},
\cite{G96}, \cite{G09}). Define
 $$
\aligned
\text{when }k\ge 0,\;\varphi_k(x_n,\sigma)&=
\cases (2\sigma)^{\frac 12}(\sigma-\partial_{x_n})^k(x_n^ke^{-x_n\sigma})/k!
&\text{ for } x_n\ge 0,\\
0 &\text{ for } x_n< 0;
\endcases\\
\text{when }k<0,\;\varphi_k(x_n,\sigma)&=\varphi_{-k-1}(-x_n,\sigma).
\endaligned
$$
Here $\sigma$ can be any positive number. In the considerations of symbols, we shall take it as
$$
\sigma (\xi ')=[\xi'],  \text{ a smooth positive function of $\xi '$ that 
equals $|\xi '|$  for $|\xi '|\ge 1$.}\tag3.1
$$ 
The functions form an orthonormal basis of $L_2({\Bbb R})$; those with
$k\ge 0$ span $L_2(\rp)$, and those with
$k< 0$ span $L_2(\rmi)$. Their
Fourier transforms are
$$
\hat\varphi _k(\xi _n,\sigma )
=(2\sigma )^{\frac12}\frac{(\sigma -i\xi _n)^k}{(\sigma
+i\xi _n)^{k+1}},\text{ where we denote }\hat\varphi '_k(\xi _n,\sigma )=\frac{(\sigma -i\xi _n)^k}{(\sigma
+i\xi _n)^{k+1}},
\tag 3.2
$$
forming orthogonal bases of 
$L_2(\Bbb R)$. 

The point is that the $\varphi _k(x_n,\sigma )$, $k\ge 0$, belong to $\Cal S_+$
and span its elements in rapidly decreasing series, in the sense that the following two systems 
of seminorms on the Fr\'echet space $\Cal S_+$ are equivalent:
$$
\aligned
&\|x_n^lD_{x_n}^{l'}u(x_n)\|_{L_2(\rp)}, \; l,l'\in{\Bbb N}_0,\text{ resp.}\\
&\|\{\ang j^Mb_j\}\|_{\ell_2({\Bbb N}_0)},\; M\in{\Bbb N}_0,\text{ where
}u(x_n)=\sum_{j\in{\Bbb N}_0}b_j\varphi _j(x_n,\sigma );
\endaligned\tag3.3
$$
cf.\ e.g.\ \cite{G96}, Lemma 2.2.1, or \cite{G09}, Lemma 10.14. This
holds in such a way that an estimate in one of the systems is
dominated by a fixed finite set of estimates in the other system. (A
sequence $\{b_j\}$ is called rapidly decreasing, when $\sup_j\ang
j^M|b_j|<\infty $ for all $M$.)

Similarly, the products $\varphi _j(x_n,\sigma )\varphi _k(y_n,\sigma
)$ belong to  $\Cal S_{++}$, 
and here there
is equivalence of the systems of seminorms 
$$
\aligned
&\|x_n^lD_{x_n}^{l'}y_n^mD_{y_n}^{m'}u(x_n,y_n)\|_{L_2(\rpp)}, \;
l,l',m,m'\in{\Bbb N}_0,\text{ resp.}\\
&\|\{\ang j^M\ang k^{M'}c_{jk}\}\|_{\ell_2({\Bbb N}_0\times{\Bbb N}_0)},\; M,M'\in{\Bbb N}_0,\text{ where
}u(x_n,y_n)=\sum_{j,k\in{\Bbb N}_0}c_{jk}\varphi _j(x_n,\sigma )\varphi _k(y_n,\sigma ),
\endaligned\tag3.4
$$
again with finite interdependence.

There are the following rules for differentiations in $\xi '$ and $\xi
_n$, when $\sigma $ is as in (3.1) (then $\partial _{\xi _j}\sigma
=\xi _j\sigma ^{-1}$ for $|\xi '|\ge 1$):
$$
\aligned
\partial_{\xi_j}\hat\varphi_k(\xi_n,\sigma)
&=\big(k\hat\varphi_{k-1}-
\hat\varphi_{k}
-(k+1)
\hat\varphi_{k+1}\big)(2\sigma)^{-1}\partial_{\xi_j}\sigma ,\quad
j<n,\\
\partial_{\xi_n}\hat\varphi_k(\xi_n,\sigma)
&=-i\big(k\hat\varphi_{k-1}+
(2k+1)\hat\varphi_{k}
+(k+1)
\hat\varphi_{k+1}\big)(2\sigma)^{-1} .
\endaligned
\tag 3.5
$$

Let us also recall that, with $e^\pm$ denoting extension by zero for
$x_n\lessgtr 0$,
Fourier transformation gives the
space $\Cal H^+=\Cal F_{x_n\to \xi _n}(e^+\Cal S_+)$; its conjugate
space is $\Cal H^-_{-1}=\Cal F_{x_n\to \xi _n}(e^-\Cal S_-)$; and both are
contained in $\Cal H=\Cal H^{+}\dot+ \Cal H^-_{-1}\dot+{\Bbb C}[t]$,
where ${\Bbb C}[t]$ consists of the polynomials on ${\Bbb R}$. The
projections of $\Cal H$ to the components $\Cal H^+$ resp.\ $\Cal H^-_{-1}$ are
denoted $h^+$ resp.\ $h^{-}_{-1}$. We refer to the indicated books
for more information; this complex notation will not be important in the
present paper.

We are now ready to give the technical definitions of the nonsmooth
symbol spaces (that may be taken lightly in a first reading). Boundary operator symbols with $C^\tau $-smoothness are defined as in Abels
\cite{A05}, and we in addition formulate the estimates in
terms of Laguerre expansions.
\proclaim{Definition 3.1} Let $d\in{\Bbb R}$ and $\tau >0$.
The space $C^\tau S^{d}_{1,0}({\Bbb R}^{n-1}\times{\Bbb R}^{n-1};\Cal S_+)$ consists
of the
functions $\tilde f(x',x_n,\xi' )$  of $x',\xi' \in{\Bbb R}^{n-1}$, $x_n\in\rp$, such that
$$
\|x_n^lD_{x_n}^{l'}D_{\xi '}^\alpha \tilde f(x',x_n,\xi ')\|_{C^\tau
({\Bbb R}^{n-1}, L_2(\rp))}\text{ is }O(\ang{\xi '}^{d+\frac12-l+l'-|\alpha
|})\text{ for  }l,l'\in{\Bbb N}_0,\alpha \in{\Bbb N}_0^{n-1}.\tag3.6 
$$

Equivalently, with  $f=\Cal F_{x_n\to\xi _n}(e^+\tilde f)$,
$$
\|h^{+}\xi _n^{l'}D_{\xi _n}^{l'}D_{\xi '}^\alpha  f(x',\xi ',\xi _n)\|_{C^\tau
({\Bbb R}^{n-1}, L_2({\Bbb R}))}\text{ is }O(\ang{\xi '}^{d+\frac12-l+l'-|\alpha
|})\text{ for  }l,l'\in{\Bbb N}_0,\alpha \in{\Bbb N}_0^{n-1}.
$$

Likewise equivalently, with $\tilde f(x',x_n,\xi ')=\sum_{k\in{\Bbb N}_0}b_k(x',\xi ')\varphi
_k(x_n,\sigma (\xi '))$,
$$
\|\{\ang k^M D_{\xi '}^\alpha b_k(x',\xi ')\}\|_{C^\tau
({\Bbb R}^{n-1}, \ell_2({\Bbb N}_0))}\text{ is }O(\ang{\xi '}^{d+\frac12-|\alpha
|})\text{ for  }M\in{\Bbb N}_0,\alpha \in{\Bbb N}_0^{n-1},\tag3.7
$$
i.e., the coefficient sequence $\{b_k(x',\xi ')\}_{k\in {\Bbb N}_0}$
is rapidly decreasing in $C^\tau S^{d+\frac12}_{1,0}({\Bbb
R}^{n-1}\times{\Bbb R}^{n-1})$.

Introduce also for $N\in{\Bbb N}_0$ the notation  $C^\tau S^{d}_{1,0}({\Bbb
R}^{n-1}\times{\Bbb R}^{n-1},N;\Cal S_+)$ for the space of functions
$\tilde f(x',x_n,\xi ')$ satisfying the estimates {\rm (3.6)} for
$l,l'\in{\Bbb N}_0$, $|\alpha |\le N$; then  \linebreak$C^\tau S^{d}_{1,0}({\Bbb
R}^{n-1}\times{\Bbb R}^{n-1};\Cal S_+)=\bigcap_{N\in{\Bbb N}_0} C^\tau S^{d}_{1,0}({\Bbb
R}^{n-1}\times{\Bbb R}^{n-1},N;\Cal S_+)$.

\endproclaim

The functions $\tilde f$ and $f$ are said be {\bf polyhomogeneous},
when there moreover is an
asymptotic expansion of $f$ in functions $f_{d-j}(x',\xi )$ homogeneous of degree
$d-j$ in $(\xi ',\xi _n)$ for $|\xi '|\ge 1$, in a similar way as in Definition
2.2.
 For the subspaces of polyhomogeneous symbols
we use the notation $C^\tau S^{d}$ instead of $C^\tau S^{d}_{1,0}$.

The functions $\tilde f$ and $f$ serve as symbol-kernels resp.\
symbols of {\bf Poisson operators of order} $d+1$, here they are
usually denoted $\tilde k$ resp.\ $ k$. First there is the definition
of an operator with respect to the $x_n$-variable:
$$
\operatorname{OPK}_n(\tilde k)v=
\tilde k(x',x_n,\xi
')\cdot v,\text{ also denoted }k(x',\xi ',D_n)v,
\tag 3.8
$$
going from ${\Bbb C}$ to $\Cal S_+$ for each $(x',\xi ')$; this is the
{\it boundary symbol operator}. Then the full operator is defined
for $v\in\Cal S({\Bbb R}^{n-1})$ by using the pseudodifferential
definition with respect to $(x',\xi ')$ (denoted $\operatorname{OP}'$):
$$
\operatorname{OPK}(\tilde k
)v=\operatorname{OP}'\operatorname{OPK}_n(\tilde k)v=\int_{{\Bbb R}^{2n-2}} e^{i(x'-y')\cdot
\xi '}\tilde k(x',x_n,\xi ') v(y')\,dy'\d\xi '.\tag3.9
$$
(One also writes $\operatorname{OPK}(\tilde k)$ as
$\operatorname{OPK}( k)$, and $\operatorname{OPK}_n(\tilde k)$ as  $\operatorname{OPK}_n( k)$.)

The same classes of functions $\tilde f$ serve as symbol-kernels 
for {\bf trace operators of order $d$ and class zero}; here
they are usually denoted $\tilde t$. The associated symbol $t$ is the
conjugate Fourier transform of $e^+\tilde t$ in $x_n$, $t(x',\xi ',\xi
_n)=\int_0^\infty e^{ix_n\xi _n}\tilde t(x',x_n,\xi ')\,
dx_n=\overline{\Cal F}_{x_n\to \xi _n}e^+\tilde t$. The definition
of an operator with respect to the $x_n$-variable is, for $u\in \Cal S_+$:
$$
\operatorname{OPT}_n(\tilde t)u=
\int_0^\infty \tilde t(x',x_n,\xi
')u(x_n)\,dx_n,\text{ also denoted }t(x',\xi ',D_n)u,\tag 3.10
$$
going from $\Cal S_+$ to ${\Bbb C}$ for each $(x',\xi ')$; the
boundary symbol operator. Then the full operator is defined
for
$u\in \Cal S(\crnp)=\Cal S({\Bbb R}^n)|_{\crnp}$ by:
$$
\operatorname{OPT}(\tilde t)u=\operatorname{OP}'\operatorname{OPT}_n(\tilde t)u=\int_{{\Bbb R}^{2n-2}}\int_0^\infty  e^{i(x'-y')\cdot
\xi '}\tilde t(x',y_n,\xi ')u(y',y_n)\,dy_ndy'\d\xi '.\tag3.11
$$

\proclaim{Definition 3.2} Let $d\in{\Bbb R}$ and $\tau >0$.
The space $C^\tau S^{d}_{1,0}({\Bbb R}^{n-1}\times{\Bbb R}^{n-1};\Cal S_{++})$ consists
of the
functions $\tilde g(x',x_n,y_n,\xi' )$  of $x',\xi' \in{\Bbb
R}^{n-1}$, $x_n, y_n\in\rp$, such that
$$
\|x_n^lD_{x_n}^{l'}y_n^mD_{y_n}^{m'}D_{\xi '}^\alpha \tilde g(x',x_n,y_n,\xi ')\|_{C^\tau
({\Bbb R}^{n-1}, L_2(\rpp))}\text{ is }O(\ang{\xi '}^{d+1-l+l'-m+m'-|\alpha
|})
\tag3.12 
$$
for  $l,l',m,m'\in{\Bbb N}_0,\alpha \in{\Bbb N}_0^{n-1}$.

Equivalently, with 
$ g=\Cal
F_{x_n\to\xi _n}\overline{\Cal F}_{y_n\to\eta
_n}(e^+_{x_n}e^+_{y_n}\tilde g)$,
$$
\|h^{+}_{\xi _n}h^-_{-1,\eta _n}\xi _n^{l'}D_{\xi _n}^{l}\eta _n^{m'}D_{\eta _n}^{m}D_{\xi '}^\alpha  g(x',\xi ',\xi _n,\eta _n)\|_{C^\tau
({\Bbb R}^{n-1}, L_2({\Bbb R}^2))}\text{ is }O(\ang{\xi '}^{d+1-l+l'-m+m'-|\alpha
|}),
$$
 for  $l,l', m, m'\in{\Bbb N}_0,\alpha \in{\Bbb N}_0^{n-1}$.

Likewise equivalenty, with $\tilde g(x',x_n,y_n,\xi ')=\sum_{j,k\in{\Bbb N}_0}c_{jk}(x',\xi ')\varphi
_j(x_n,\sigma )\varphi _k(y_n,\sigma )$,
$$
\|\{\ang j^M\ang k^{M'} D_{\xi '}^\alpha c_{jk}(x',\xi ')\}\|_{C^\tau
({\Bbb R}^{n-1}, \ell_2({\Bbb N}_0\times {\Bbb N}_0))}\text{ is }O(\ang{\xi '}^{d+1-|\alpha
|}),\tag3.13
$$
for  $M,M'\in{\Bbb N}_0,\alpha \in{\Bbb N}_0^{n-1}$, i.e., the coefficient sequence $\{c_{jk}(x',\xi ')\}_{j,k\in {\Bbb N}_0}$
is rapidly decreasing in $C^\tau S^{d+1}_{1,0}({\Bbb
R}^{n-1}\times{\Bbb R}^{n-1})$.

Introduce also the notation $C^\tau S^{d}_{1,0}({\Bbb R}^{n-1}\times{\Bbb
R}^{n-1},N;\Cal S_{++})$ for the space of functions
$\tilde g(x',x_n,y_n,\xi' )$ satisfying {\rm (3.12)} for
$l,l',m,m'\in{\Bbb N}_0$, $|\alpha |\le N$.
\endproclaim

The functions $\tilde g$ and $g$ are said be {\bf polyhomogeneous},
when there moreover is an
asymptotic expansion of $g$ in functions $g_{d-j}(x',\xi ',\xi _n,\eta _n)$ homogeneous of degree
$d-j$ in $(\xi ',\xi _n,\eta _n)$ for $|\xi '|\ge 1$, in a similar way
as in Definition 2.2.

The functions $\tilde g$ and $g$ serve as symbol-kernels resp.\
symbols of {\bf singular Green operators of order $d+1$ and class
zero}.  The definition
of an operator with respect to the $x_n$-variable is
$$
\operatorname{OPG}_n(\tilde g)u=
\int_0^\infty \tilde g(x',x_n,y_n,\xi
')u(y_n)\,dy_n,\text{ also denoted }g(x',\xi ',D_n)u,\tag 3.14
$$
acting in $\Cal S_+$ for each $(x',\xi ')$; the
boundary symbol operator. Then the full operator is defined
for
$u\in \Cal S(\crnp)$ by:
$$
\operatorname{OPG}(\tilde g)u=\operatorname{OP}'\operatorname{OPG}_n(\tilde g)u
=\int _{{\Bbb R}^{2n-2}}e^{i(x'-y')\cdot
\xi '}\int_0^\infty \tilde g(x',x_n,y_n,\xi ') u(y',y_n)\,dy_n dy'\d\xi '.\tag3.15
$$

The operators defined as in (3.9), (3.15) can for precision be 
said to be
{\it in $x'$-form}, to distinguish this from the case where the
functions $\tilde k$, $\tilde t$, $\tilde g$ in the integrals depend on $y'$ in
the place of $x'$.  Such cases also define Poisson, trace and singular
Green operators; the operators are said to be {\it in $y'$-form},
denoted $\operatorname{OPK}(\tilde k(y',x_n,\xi '))$, etc. Also
$(x',y')$-forms can occur. The
adjoint of the Poisson operator $\operatorname{OPK}(\tilde
f(x',x_n,\xi '))$ is the trace operator $\operatorname{OPT}(\overline{\tilde f}(y',x_n,\xi '))$.

In the following we often leave out ${\Bbb R}^{n-1}\times{\Bbb
R}^{n-1}$ from the notation for the symbol spaces.

For the consideration of the full operators as operator-valued
$\psi $do's, the following properties of the boundary symbol operators
will be very useful:

\proclaim{Lemma 3.3} Let $\tilde k\in C^\tau S^{d-1}_{1,0}({\Bbb
R}^{n-1}\times{\Bbb R}^{n-1};\Cal S_+)$, $\tilde t\in C^\tau S^{d}_{1,0}({\Bbb
R}^{n-1}\times{\Bbb R}^{n-1};\Cal S_+)$, and   $\tilde g\in C^\tau S^{d-1}_{1,0}({\Bbb
R}^{n-1}\times{\Bbb R}^{n-1};\Cal S_{++})$. Then the boundary
symbol operators  $k(x',\xi ',D_n) $, $t(x',\xi ',D_n)$ and
$g(x',\xi ',D_n)$ (cf.\ {\rm (3.8), (3.10), (3.14)})
satisfy, for $s',s''\ge 0$:
$$
\aligned
k(x',\xi ',D_n)&\in C^\tau S^{d-\frac12+s'}_{1,0}({\Bbb R}^{n-1}\times{\Bbb
R}^{n-1}; \Cal L({\Bbb C}, H^{s'}(\rp)),\\
t(x',\xi ',D_n)&\in C^\tau S^{d+\frac12+s'}_{1,0}({\Bbb R}^{n-1}\times{\Bbb
R}^{n-1}; \Cal L( H^{-s'}_0(\crp),{\Bbb C})),
\\
g(x',\xi ',D_n)&\in C^\tau S^{d+s'+s''}_{1,0}({\Bbb R}^{n-1}\times{\Bbb
R}^{n-1}; \Cal L( H^{-s'}_0(\crp),H^{s''}(\rp)).
\endaligned\tag3.16
$$
 Here, when $M\in{\Bbb N}_0$, the Poisson boundary symbol operator  seminorms 
$$
\sup_{\xi '}\ang{\xi '}^{-d+\frac12+|\alpha |-s'}\|D_{\xi '} ^\alpha k(x',\xi ',D_n)\|_{C^\tau ({\Bbb
R}^{n-1},\Cal L({\Bbb C}, H^{s'}({\Bbb R}_+)))},\quad |\alpha |\le N,
s'\le M,\tag3.17
$$ 
are dominated by the system of symbol seminorms
$$
\sup_{\xi '}\ang{\xi '}^{-d+\frac12+|\alpha |-l}\|D_{x_n}^{l}D_{\xi '} ^\alpha \tilde k(x',x_n,\xi ')\|_{C^\tau ({\Bbb
R}^{n-1},L_2(\rp))},\quad |\alpha |\le N,l\le M
.\tag3.18
$$
They are also dominated by the system of  seminorms in terms of Laguerre expansions 
$\tilde k(x',x_n,\xi ')=\sum_{j\in{\Bbb N}_0}b_j(x',\xi ')\varphi
_j(x_n,\sigma (\xi '))$: $$
\sup_{\xi '}\ang{\xi '}^{-d+\frac12+|\alpha
|}\|\{\ang j^{M+1}  D_{\xi '}^\alpha b_j(x',\xi ')\}_{j\in{\Bbb N}_0}\|_{C^\tau
({\Bbb R}^{n-1}, \ell_2({\Bbb N}_0))},\quad 
|\alpha |\le N.\tag3.19
$$ 

Likewise, the trace and singular Green boundary symbol operator
seminorms  
$$
\aligned
&\sup_{\xi '}\ang{\xi '}^{-d-\frac12+|\alpha |-s'}\|D_{\xi '} ^\alpha t(x',\xi ',D_n)\|_{C^\tau ({\Bbb
R}^{n-1},\Cal L(H^{-s'}_0(\crp),{\Bbb C}))},\quad |\alpha |\le N,
s'\le M,\text{ resp.}\\
&\sup_{\xi '}\ang{\xi '}^{-d+|\alpha |-s'-s''}\|D_{\xi '} ^\alpha g(x',\xi ',D_n)\|_{C^\tau ({\Bbb
R}^{n-1},\Cal L(H^{-s'}_0{\crp}, H^{s''}({\Bbb R}_+)))},\quad |\alpha |\le N,
s'\le M,
\endaligned
\tag3.20
$$ 
are dominated by the systems of symbol seminorms
$$
\aligned
&\sup_{\xi '}\ang{\xi '}^{-d-\frac12+|\alpha |-l}\|D_{x_n}^{l}D_{\xi '} ^\alpha \tilde t(x',x_n,\xi ')\|_{C^\tau ({\Bbb
R}^{n-1},L_2(\rp))},\; |\alpha |\le N,l\le M,\text{ or}\\
&\sup_{\xi '}\ang{\xi '}^{-d-\frac12+|\alpha
|}\|\{\ang j^{M+1}  D_{\xi '}^\alpha b_j(x',\xi ')\}_{j\in{\Bbb N}_0}\|_{C^\tau
({\Bbb R}^{n-1}, \ell_2({\Bbb N}_0))},\; 
|\alpha |\le N, \text{ resp.}\\
&\sup_{\xi '}\ang{\xi '}^{-d+|\alpha |-l}\|D_{x_n}^{l}D_{y_n}^{m}D_{\xi '} ^\alpha \tilde g(x',x_n,y_n,\xi ')\|_{C^\tau ({\Bbb
R}^{n-1},L_2(\rpp))},\; |\alpha |\le N,l,m\le M,\text{ or}\\
&\sup_{\xi '}\ang{\xi '}^{-d+|\alpha
|}\|\{\ang {(j,k)}^{2M+1}  D_{\xi '}^\alpha c_{jk}(x',\xi ')\}_{j,k\in{\Bbb N}_0}\|_{C^\tau
({\Bbb R}^{n-1}, \ell_2({\Bbb N}^2_0))},\; 
|\alpha |\le N.\\
\endaligned
\tag3.21$$ 
Here $\tilde t=\sum_{j\in{\Bbb N}_0}b_j\varphi
_j(x_n,\sigma )$ and $\tilde g=\sum_{j,k\in{\Bbb N}_0}c_{jk}\varphi
_j(x_n,\sigma )\varphi _k(y_n,\sigma )$. 
\endproclaim

\demo{Proof} 
This is a variant of a result shown (in more generality including $L_p$-spaces
and weighted norms) in Abels \cite{A05}, Lemma
4.6, as a generalization of \cite{GK93}, Th.\ 3.7. 

The present $L_2$-version is straightforward to show. Consider $\tilde k$. 
The estimates
$$
\ang{\xi '}^{-d+\frac12-l}\|D_{x_n}^{l}\tilde k(x',x_n,\xi
')\|_{C^\tau ({\Bbb R}^{n-1}, L_2(\rp))}=O(1),\quad l\le M,
$$
imply 
$$
\ang{\xi '}^{-d+\frac12-M}\|\tilde k(x',x_n,\xi
')\|_{C^\tau ({\Bbb R}^{n-1}, H^M(\rp))}=O(1),
$$
so since $k(x',\xi ',D_n)$ is multiplication by $\tilde k$, the
finiteness of (3.17)
follows for $\alpha =0$, $s'=M$. A similar treatment of $D_{\xi
'}^\alpha \tilde k$ includes general $\alpha $. Noninteger values
$s'\in [0,M]$ are included by interpolation, since $H^{s'}$
interpolates between $H^0$ and $H^M$; here the symbol seminorms for
$l\le M$ suffice for the estimates. This shows the statement on
domination of the seminorms (3.17) by those in (3.18).

The second statement on domination follows from the fact that by
\cite{G96}, (2.2.20),
$$
\ang{\xi '}^{-d+\frac12-l}\|D_{x_n}^{l}\tilde k(x',x_{n},\xi
')\|_{L_2(\rp)}\le c\|\{\ang j^{l+\varepsilon }b_j(x',\xi ')\}_{j\in{\Bbb N}_0}\|_{\ell_2({\Bbb N}_0)}.
$$

The proofs for the other types of operators are similar (for $\tilde
t$ one can simply note that $t(x',\xi ',D_n)$ is the adjoint of
the Poisson operator $v\mapsto \overline{\tilde t}(x',x_n,\xi ')v$, and
$H^{-s'}_0(\crp)$ is the dual space of $H^{s'}(\rp)$).
\qed
\enddemo

One of the consequences derived in \cite{A05} is that
the operators of order $d$ have the Sobolev space continuities:$$\aligned
\operatorname{OPK}(\tilde k)&\colon H^{s+d-\frac12}({\Bbb R}^{n-1})\to
H^s(\rnp)\text{ for }|s|<\tau ,\\
\operatorname{OPT}(\tilde t)&\colon H^{s+d}(\rnp)\to
H^{s-\frac12}({\Bbb R}^{n-1})
\text{ for }|s-\tfrac12|<\tau ,\,  s+d>-\tfrac12,\\
\operatorname{OPG}(\tilde g)&\colon H^{s+d}(\rnp)\to
H^{s}(\rnp)
\text{ for }|s|<\tau , \, s+d>-\tfrac12;
\endaligned\tag3.22$$
when $\tilde k$, $\tilde t$ and $\tilde g$ are as in Lemma 3.3.
We return to the proof below in Theorem 3.8. 

\example{Remark 3.4} The special case where the operators map from $H^s$-space to
$H^s$-space for $|s|<\tau$, is particularly convenient in composition rules:
$$\aligned
\operatorname{OPK}(\tilde k)&\colon H^{s}({\Bbb R}^{n-1})\to
H^s(\rnp)\text{ for }|s|<\tau ,\text{ when }d=\tfrac12,\\
\operatorname{OPT}(\tilde t)&\colon H^{s}(\rnp)\to
H^{s}({\Bbb R}^{n-1})
\text{ for }|s|<\tau ,\,  s>-\tfrac12,\text{ when }d=-\tfrac12,\\
\operatorname{OPG}(\tilde g)&\colon H^{s}(\rnp)\to
H^{s}(\rnp)
\text{ for }|s|<\tau , \, s>-\tfrac12,\text{ when }d=0;
\endaligned\tag3.23$$
here $\tilde k$ and $\tilde t\in C^\tau
S^{-\frac12}_{1,0}(\Cal S_+)$, $\tilde g\in C^\tau S^{-1}_{1,0}(\Cal
S_{++})$. (For the statement on $\operatorname{OPT}(\tilde t)$ we
replaced $s-\frac12$ in (3.22) by $s$.) We say that these operators are {\it of neutral order}.

The order conventions, introduced originally by Boutet de Monvel in
\cite{B71}, may seem a little confusing. They respect the principle
that a composition of two operators of orders $d_1$ resp.\ $d_2$ is of order
$d_1+d_2$, but the order is not preserved when one passes from a Poisson
operator of order $d$ to its adjoint, which is a  trace operator of order $d-1$
and class $0$. To compensate for this phenomenon, one might think of
redefining the orders by adding or subtracting $\frac12$, but this is not really helpful,
 since in considerations of mappings
between $L_p$-based Sobolev spaces, the role of $\frac12$ is taken over
by $\frac1p$ and $\frac1{p'}$.
\endexample

\subhead 3.2 Rules of calculus \endsubhead

Let us now recall some composition rules. We here let
$a(x',\xi ',D_n)$ play the role of any of the boundary symbol
operators introduced above, such that the resulting full operator can
be written as  $A=\operatorname{OP}'a(x',\xi ',D_n)$. We also include
functions $s(x',\xi ')$ that are symbols of  $\psi $do's on ${\Bbb
R}^{n-1}$, and on the boundary symbol level simply act as multiplications. 
The composition of two boundary symbol operators $a_1$ and $a_2$ is
denoted $a_1\circ_n a_2$.
The compositions can of course only be applied when $A_1$ is
defined on the range of $A_2$ (which can be a Sobolev space over ${\Bbb R}^{n-1} $
or $\rnp$). The notation $\tilde a_1\circ_n \tilde a_2=\tilde a_3$ is
sometimes also used with
the $\tilde a_i$ denoting the corresponding symbol-kernels. Full
details are found in  \cite{G96, G09} or \cite{A05}; some examples are:
$$\aligned
\tilde k\circ_n\tilde t&=\tilde k(x',x_n,\xi ')\tilde t(x',y_n,\xi
'),\\
\tilde t\circ_n \tilde k&=\int _0^\infty \tilde t(x',x_n,\xi ')\tilde k(x',x_n,\xi
')\,dx_n,\\
\tilde g\circ_n \tilde k&=\int _0^\infty \tilde g(x',x_n,y_n,\xi ')\tilde k(x',y_n,\xi
')\,dy_n,\\
\tilde t\circ_n \tilde g&=\int _0^\infty \tilde t(x',x_n,\xi ')\tilde g(x',x_n,y_n,\xi
')\,dx_n,
\endaligned\tag3.24$$
$$
\aligned
\tilde g\circ_n \tilde g'&=\int _0^\infty \tilde g(x',x_n,z_n,\xi ')\tilde g'(x',z_n,y_n,\xi
')\,dz_n,\\
\tilde k\circ_n s&=\tilde k(x',x_n,\xi ')s(x',\xi '),\\
s\circ_n\tilde t&=s(x',\xi ')\tilde t(x',x_n,\xi
').
\endaligned$$

As usual, $\operatorname{OP}'(a_1\circ a_1)$ is a good
approximation to $A_1A_2$; the following theorem gives more
information on this for the H\"o{}lder-smooth symbol classes.

\proclaim{Theorem 3.5}  
For $i=1,2$, let $\tilde k_i(x',x_n,\xi
')$ be Poisson symbol-kernels of order $d_j+\frac12$, let $\tilde
t_i(x',x_n,\xi ')$ be  trace symbol-kernels
 of order $d_i-\frac12$ and class $0$, let 
 $\tilde g_i(x',x_n,y_n,\xi
')$ be singular Green symbol-kernels of order $d_i$ and class $0$,
and let $s_i(x',\xi ')$ be
$\psi $do symbols on ${\Bbb R}^{n-1}$ of orders $d_i$,
with  $\tau _i$-smoothness, respectively, defining operators $K_i$,
$T_i$, $G_i$, $S_i$. Let $d=d_1+d_2$ and $\tau =\min\{\tau
_1,\tau _2\}$, and let $\theta \in \,]0,1 [\,$ with $\theta <\tau $. Then the composed operators
satisfy$$\aligned
A_1A_2&-\operatorname{OP}'(a_1\circ_n
a_2)\colon H^{s+d-\theta }\to
H^{s},\\
\text{ when } &|s|<\tau ,\; s>\theta -\tau _2,\; \theta -\tau
_2<s+d_1<\tau _2;
\endaligned\tag3.25
$$
in addition the class condition $s+d-\theta >-\frac12$ must be
satisfied if  $A_2$ is a trace or singular Green
operator. (The $H^t$ stand for Sobolev spaces over ${\Bbb R}^{n-1}$ or
$\rnp$ depending on the context where they are used.) 

\endproclaim

The rule in the case where $a_1$ and $a_2$ are $\psi $do's on ${\Bbb
R}^{n-1}$ is known from Taylor \cite{T91},
Prop.\ I 2.1D. The other rules (on $\psi $dbo's) are
proved in Abels \cite{A05}, Th.\ 4.13, by use of an extension of
Taylor's result to vector-valued operators, cf.\ also
\cite{A05a}. The above statement differs from the formulation in
\cite{A05} by referring to orders
$d_i\pm\frac12$ for Poisson and trace operators; this allows a unified
formulation. One can also
describe rules where a trace operator or singular Green operator
contains standard traces $\gamma _j$ (i.e.,  has
positive class); they can be deduced from the above by combination
with mapping properties of the $\gamma _j$.

We moreover need rules for compositions with $\psi $do's $P$ on ${\Bbb
R}^n$ satisfying the transmission condition at $x_n=0$, and truncated
to $\rnp$ as $P_+=r^+Pe^+$ (where $r^\pm$ denotes restriction to
$\rnpm$, $e^\pm$ stands for extension by zero on ${\Bbb R}^n\setminus\rnpm$). Here the
composition rules are more complicated, already on the symbol level. 
Moreover, the symbol $p(x,\xi
)$ of $P$ is in general $x_n$-dependent, whereas the boundary symbol
operator $p(x',\xi ',D_n)_+$ is defined from the symbol at $x_n=0$. When
one of the factors $A_i$ in a composition is of the 
form $P_+$, 
it is 
the true operator that enters as $A_i$ and the $x_n$-independent
boundary symbol operator that enters as $a_i$. It is proved in Section
5 of \cite{A05} that the statement in Theorem 3.5
is valid also for these cases:

\proclaim{Theorem 3.6} The conclusion of Theorem {\rm 3.5} holds also
when $G_1$ or $G_2$ is replaced by a $C^{\tau _i}$-smooth  $\psi $do $P_{i,+}$ of
order $d_i$ satisfying the transmission condition at $x_n=0$. 
\endproclaim
 
We here view the truncated $\psi $do's as being of class 0. There is a
refinement allowing operators of negative class, but we shall not need
it here and refer to the quoted works for details on this.

For a composition of two truncated $\psi $do's $P_{1,+}$ and $P_{2,+}$
one uses the formula
$$
P_{1,+}P_{2,+}=(P_{1}P_2)_+-L(P_1,P_2),\tag3.26
$$
where the composition $P_1P_2$ follows a rule as for $S_1S_2$ in
Theorem 3.5, and the singular Green term $L(P_1,P_2)$ is treated by use
of other rules in Theorem 3.5, supplied with considerations of
standard trace operators $\gamma _j$. The outcome is essentially as in (3.25)
 (more details are given in \cite{A05},
Sect.\ 5.3).

For the purpose of spectral estimates we need a sharpening of the information in (3.22), (3.23), taking
into account how many symbol seminorm estimates are needed for the operator
norm estimates. This uses the following generalization of Marschall's
result Theorem 2.3:

\proclaim{Theorem 3.7} Let $H_0$ and $H_1$ be Hilbert spaces, and
consider operator-valued symbols $p\in C^\tau S^d_{1,0}({\Bbb
R}^m\times{\Bbb R}^{m}, N; \Cal L(H_0,H_1))$, where the symbol
seminorms {\rm(2.8)} are taken in $C^\tau ({\Bbb R}^m, \Cal
L(H_0,H_1))$. For $N\ge m+1$, $\operatorname{OP}(p)$ is continuous
$$
\operatorname{OP}(p)\colon H^{s+d}({\Bbb R}^m,H_0)\to H^{s}({\Bbb
R}^m,H_1)\text{ for }|s|<\tau .
$$
 \endproclaim

The proof is given in the Appendix (written by
Helmut Abels). We use the theorem to show that a fixed finite set of
symbol seminorms suffices for estimates of operator norms.

\proclaim{Theorem 3.8} Let $\tau \in \,]0,1[\,$. 
Let $\tilde k$, $\tilde t$ and $\tilde g$ be as in Lemma {\rm 3.3}.

$1^\circ$ For each $|s|<\tau $, the norm of $\operatorname{OPK}(\tilde
k)$ in {\rm (3.22)} is bounded by a finite system of seminorms:  
$$
\aligned
&\|\operatorname{OPK}(\tilde k)\|_{\Cal L(H^{s+d-\frac12}({\Bbb R}^{n-1}), H^s(\rnp))}
\\
&\le C_s\sup_{\xi '}\ang{\xi '}^{|\alpha |-d+\frac12}\|D_{x_n}^lD_{\xi
'}^{\alpha }\tilde k\|_{C^\tau ({\Bbb R}^{n-1},L_2(\rp))},\;l=0,1,\; |\alpha
|\le n.
\endaligned\tag3.27
$$
The symbol seminorms in {\rm (3.27)} are
estimated by the following Laguerre symbol seminorms, where $\tilde k(x',x_n,\xi ')=\sum_{k\in{\Bbb N}_0}b_k(x',\xi ')\varphi _k(x_n,\sigma (\xi '))$,
$$
\sup_{\xi '}\ang{\xi '}^{|\alpha |-d+\frac12}\|\{\ang k^2D_{\xi
'}^\alpha b_{k}\}\|_
{C^r({\Bbb R}^{n-1},\ell_2({\Bbb N}_0))},\, |\alpha |\le n.\tag3.28
$$

$2^\circ$ 
 For each $|s'-\frac12|<\tau $, resp.\ $|s|<\tau $, the norms of
$\operatorname{OPT}(\tilde t)$  and $\operatorname{OPG}(\tilde g)$ 
in {\rm (3.22)} are bounded by  finite systems of seminorms:  
$$
\aligned
&\|\operatorname{OPT}(\tilde t)\|_{\Cal L(H^{s'+d}(\rnp), H^{s'-\frac12}(\rnp))}
\\
&\le C_{s'}\sup_{\xi '}\ang{\xi '}^{|\alpha |-d-\frac12}\|D_{x_n}^lD_{\xi
'}^{\alpha }\tilde t\|_{C^\tau ({\Bbb R}^{n-1},L_2(\rp))},\;l=0,1,\, |\alpha
|\le n,\\
&\|\operatorname{OPG}(\tilde g)\|_{\Cal L(H^{s+d}(\rnp), H^s(\rnp))}
\\
&\le C_s\sup_{\xi '}\ang{\xi '}^{|\alpha |-d}\|D_{x_n}^lD_{y_n}^{m}D_{\xi
'}^{\alpha }\tilde g\|_{C^\tau ({\Bbb R}^{n-1},L_2(\rpp))},\;l, m=0,1,\, |\alpha
|\le n.
\endaligned\tag3.29
$$
The operator norms can also be estimated in terms of finite systems of
Laguerre seminorms, as in Lemma {\rm 3.3}.
\endproclaim

\demo{Proof}  The proof follows that of
\cite{A05}, Th.\ 4.8 in a simplified version, but taking the dependence on
specific finite seminorm systems into account.

Consider $K=\operatorname{OPK}(\tilde k)=\operatorname{OP}'(k(x',\xi
',D_n))$. By Lemma 3.3, $k(x',\xi ',D_n)$ is in  \linebreak$C^\tau S^{d-\frac12+s'}_{1,0}({\Bbb R}^{n-1}\times{\Bbb
R}^{n-1}; \Cal L({\Bbb C}, H^{s'}(\rp))$ for all $s'\ge 0$. Applying
Theorem 3.7 with $N=(n-1)+1=n$, we see that $K$ is continuous 
$$
K\colon H^{s+d-\frac12+s'}({\Bbb R}^{n-1})\to H^{s}({\Bbb R}^{n-1}, H^{s'}(\rp)) \tag3.30
$$
 for $|s|<\tau $, $s'\ge 0$. Here we let $s'\in [0,1]$, then the
 operator norm is estimated by (3.18) or (3.19) with $N=n$ and
 $M=1$. If $s\in [0,\tau [\,$, observe that
$$
H^s(\rnp)=H^s({\Bbb R}^{n-1}, H^0(\rp))\cap H^0({\Bbb R}^{n-1},
H^s(\rp)) \text{ for }s\ge 0,\tag3.31
$$
and apply this to (3.30) with $(s,s')$ of the form $(s,0)$ and of the
form
$(0,s)$; then we find that 
$$
K\colon H^{s+d-\frac12}({\Bbb R}^{n-1})\to H^{s}(\crp) \tag3.32
$$
holds for $s\ge 0$. If $s\in \,]-\tau ,0]$, we obtain (3.32) by taking $s'=0$ in
(3.30) and using that $$
H^s({\Bbb R}^{n-1}, H^0(\rp))\hookrightarrow
H^s(\rnp)\text{ if }s\le 0.\tag3.33$$

Since $s'\in [0,1]$, the operator norm is estimated as asserted in the
theorem.

The proofs for $T$ and $G$ are adapted from \cite{A05} Th.\ 4.8 in a
similar way.
\qed 
\enddemo

Also
rules in Theorems 3.4 and 3.5 may be sharpened to a dependence on
only finitely many symbol seminorms, by an extension of results of
Marschall \cite{M87} to vector-valued operators, but since we do not need them in the present
paper, we  shall not pursue this here.

\example {Remark 3.9} Let us mention one more rule, namely that the
difference between operators defined in
$x'$-form resp.\ $y'$-form from the same symbol is of lower
order; it is not included in \cite{A05, A05a}. 
Since the result is not essential for the
present work, we only indicate some ingredients in a proof:
A special case of
Corollary 4.6 of \cite{M87} is: If 
$a(x,\xi )\in C^\tau S^d_{1,0}({\Bbb R}^m\times{\Bbb R}^m, N;\Cal
L(H_0,H_1))$, with $H_0=H_1={\Bbb C}$, then
$$
\operatorname{OP}(a(x,\xi ))-\operatorname{OP}(a(y,\xi ))\colon
H^{s+d-\theta }({\Bbb R}^m,H_0)\to H^s({\Bbb R}^m,H_1)\tag 3.34
$$
is bounded, provided that $
N>\tfrac32 m+2$, $\theta \in [0,1]$, $\tau >\theta $, $|s|<\tau $ and
$|s+d-\theta |<\tau $. (This allows $d$ and $s+d$ to run in a small interval
around 0.) 

One needs a generalization of the statement to general Hilbert spaces
$H_0,H_1$, and then one can prove by use of Lemma 3.3 that when
$K=\operatorname{OPK}({\tilde k}(x',x_n,\xi '))$ and
$K'=\operatorname{OPK}({\tilde k}(y',x_n,\xi '))$ are Poisson
operators 
defined in $x'$-form resp.\ $y'$-form from a symbol $\tilde
k(x',x_n,\xi ')\in C^\tau S^{-\frac12}_{1,0}({\Bbb R}^{n-1}\times
{\Bbb R}^{n-1};\Cal S_+)$, $\tau <1$, then
$$
K-K'\colon H^{s-\theta }({\Bbb R}^{n-1})\to H^s(\rnp)\text{ for
}-\tau +\theta <s<\tau . \tag3.35
$$
There are similar statements for trace and singular Green operators.

\endexample

\subhead 3.3 Laguerre boundary operators \endsubhead

Let us introduce a notation for the special Poisson and trace operators
on $\rnp$ --- Laguerre boundary operators --- defined from Laguerre functions with
$k\ge 0$ (as in
Grubb and Schrohe \cite{GS04}):
$$
\Phi _k=\operatorname{OPK}({\varphi}_k(x_n,\sigma(\xi '))),\quad \Phi ^*
_k=\operatorname{OPT}({\varphi}_k(x_n,\sigma(\xi ')));\tag 3.36
$$
here $\Phi _k$ is a Poisson operator of order $\frac12$, and the adjoint 
$\Phi _k^*$ is a trace operator of order $-\frac12$ and class 0 
(note that they are of neutral order as defined in Remark 3.4). The operator
$\Phi _k$ maps $L_2(\Bbb R^{n-1})$ continuously, in fact
isometrically, into $L_2(\rnp)$.
Moreover, because of
the orthonormality and completeness of the $\varphi _k$,
$$
\Phi ^*_j\Phi
_k=\delta _{jk}I_{{\Bbb R}^{n-1}},\quad \sum_{k\in {\Bbb N}_0}\Phi _k\Phi _k^*=I_{\rnp},\tag 3.37
$$ where $I_{{\Bbb R}^{n-1}}$ resp.\ $I_{\rnp}$ stands for the identity operator on functions on
$\Bbb R^{n-1}$ resp.\ $\rnp$. Denoting by $H_j$ the range of $\Phi _j$
in $L_2(\rnp)$, we have that the $H_j$ are mutually orthogonal closed
subspaces, and that
$$
L_2(\rnp)=\bigoplus _{j\in{\Bbb N}_0}H_j;\text{ with isometries }\Phi
_j\colon L_2({\Bbb R}^{n-1})\simto H_j,\quad \Phi _j^*\colon H_j\simto
L_2({\Bbb R}^{n-1}),\tag3.38
$$
the latter acting as an inverse of $\Phi _j$.

\proclaim{Lemma 3.10}
When $G$ is a singular Green operator on $\rnp$ of order $d$ and class
$0$, defined from a symbol-kernel $\tilde g(x',\xi ',x_n,y _n)\in C^\tau S^{d-1}_{1,0} ({\Bbb
R}^n_+\times{\Bbb R}^{n-1};\Cal S_{++})$, 
it can be written in the form
$$
G=\sum_{k\in\Bbb N_0} K_k
\Phi ^*_k
,\quad K_{k}=\operatorname{OPK}(\tilde k_{k}(x',x_n,\xi')) , \tag 3.39$$
with Poisson operators $K_k$ of order $d+\frac12$, their
symbol-kernels $(\tilde k_{k})_{k\in\Bbb N_0}$ forming a rapidly
decreasing sequence in
$C^\tau S^{d-\frac12} _{1,0}({\Bbb R}^{n-1}\times {\Bbb
R}^{n-1};\Cal S_+)$. 
Here
$$
\tilde k_{k}(x',x_n,\xi ')=\int_0^\infty \tilde g(x',x_n,y_n,\xi ')\varphi _k(y _n,\sigma (\xi
'))\,d y _n, \quad K_k=G\Phi _k.\tag3.40
$$

\endproclaim

\demo{Proof} The symbol-kernel $\tilde g$ of $G$ has a Laguerre expansion:
$$
g(x',x_n,y_n,\xi') = \sum_{j,k\in\Bbb N_0}
c_{jk}(x',\xi'){\varphi}_j(x_n,\sigma(\xi '))
\varphi_k(y_n,\sigma(\xi ')),
\tag 3.41$$
with $(c_{jk})_{j,k\in\Bbb N_0}$ rapidly decreasing in $C^\tau S^{d}_{1,0}$,
i.e., the relevant symbol seminorms on \linebreak
$\ang k^{M}\ang
j^{M'}c_{jk}$ are
bounded in $j,k$ for any $M,M'\in\Bbb N_0$. Set
$$
\tilde k_k(x',x_n,\xi ')=\sum_{j\in\Bbb N_0}
c_{jk}(x',\xi')\varphi_j(x_n,\sigma).
$$
It is a rapidly decreasing sequence of symbols in $C^\tau S^{d-\frac12} _{1,0}({\Bbb R}^{n-1}\times {\Bbb
R}^{n-1};\Cal S_+)$, in view of the estimates of the $c_{jk}$ and the
($x'$-independent) $\varphi _j$.  
Then
$$
G=\sum_{k\in\Bbb N_0} \operatorname{OPK}\bigl(\tilde k_k(x',x_n,\xi ')\bigr)\circ
\operatorname{OPT}\bigl(\varphi _k(x _n,\sigma )\bigr)=\sum_{k\in\Bbb
N_0} 
K_k\Phi _k^*,
$$
where $K_k=\operatorname{OPK} (\tilde k_k)$.
The Poisson operators $K_k$ and their symbols satisfy (3.40) in view
of the orthonormality of the $\varphi _k$ and (3.37).

This shows the assertion.\qed
\enddemo

\head 4. Spectral estimates for nonsmooth singular Green operators \endhead

\subhead 4.1 Upper estimates \endsubhead

We now turn to eigenvalue estimates for s.g.o.s. First we establish upper spectral estimates. (Both here and in the proof of
Theorem 4.3 below, the ideas from Ch.\ 4 of \cite{G84}
are used wherever convenient.)
Let $\psi $ and $\psi _1$ be functions in $C_0^\infty
({\Bbb R}^n)$, supported in a ball $B_R$. Denote
$B_R\cap\crnp=B_{R,+}$, $B_R\cap( {\Bbb R}^{n-1}\times \{0\})=B'_R$.

\proclaim{Theorem 4.1} Let $G$ be as in Lemma {\rm 3.10}, of negative order
$d=-t$.  Then $\psi  G\psi _1\in \frak S_{(n-1)/t,\infty }$. Moreover, when $G$ is
written as a series $\sum_{k\in{\Bbb N}_0}K_k\Phi _k^*$, 
$$
{\bold N}_{(n-1)/t}(\psi G\psi _1)\le C \sum_{k\in {\Bbb
N}_0}k^\gamma {\bold N}_{(n-1)/t}(\psi K_k\Phi _k^*
\psi _1)<\infty ,\tag4.1
$$
where 
$\gamma =0$ if $t<n-1$, $\gamma >t/(n-1)-1$ if $t\ge n-1$.
\endproclaim

\demo{Proof} The crucial step is to show that each $\psi K_k\Phi _k^*\psi _1$ is
compact as an operator in $L_2(B_{R,+})$, with
the desired estimate of $s$-numbers. To see this, note that  $K_k$ is
bounded from  $H^{-t}({\Bbb R}^{n-1})$ to $L_2(\rnp)$, where $\psi 
K_k$ maps into functions supported in $B_{R,+}$, with a
bound depending on the symbol estimates and $\psi $. Then the
adjoint $K_k^*\psi $ is bounded from $L_2(B_{R,+})$ to 
$H^t({\Bbb R}^{n-1})$. But we cannot conclude from this that  $\psi
K_k $ and $K_k^*\psi $ would be in $\frak S_{(n-1)/t,\infty }$, 
since there is no support restriction in ${\Bbb R}^{n-1}$.

Instead we take $\psi _1$ into the picture, and 
moreover insert a cutoff function $\zeta\in C_0^\infty ({\Bbb R}^{n-1})$ that is 1 on $B'_R$, to the right of $K_k$. Write
$$
\psi K_k\Phi ^*_k\psi _1=\psi K_k\zeta\Phi _k^* \psi _1+\psi K_k(1-\zeta)\Phi _k^*\psi _1. \tag4.2
$$
The first term will have the $\frak S_{(n-1)/t,\infty }$-property, since
$\psi K_k\zeta$ has it and the
$\Phi _k^*$ are bounded from $L_2(\rnp)$ to $L_2({\Bbb R}^{n-1})$ with
norm 1. For the second term, we focus instead on \linebreak$(1-\zeta)\Phi
^*_k\psi _1$. As shown in \cite{G84}, proof of Prop.\ 4.7, it is of order $-\infty
$ in such a way that for its estimates as a mapping of a fixed low order, the
symbol seminorms that enter have a polynomial growth in $k$ (in view
of the formulas (3.5)). Then for a fixed large $a$, the norm of
$(1-\zeta)\Phi ^*_k\psi _1$ from $H^{-a}_0(\crnp)$ to $H^{-t}({\Bbb
R}^{n-1})$ grows at most polynomially with $k$. Let $\widetilde B$ be
a smooth bounded open subset of $\rnp$ containing $B_{R,+}$. The operator $\psi
K_k(1-\zeta )\Phi _k^*\psi _1$ is bounded from $H^{-a}_0(\overline
{\widetilde B})$ to
$L_2(\widetilde B)$, 
and its adjoint is bounded from $L_2(\widetilde B )$  to
$H^{a}(\widetilde B)$, hence
belongs to $\frak S_{n/a,\infty }$. We can take $a$ so large that $n/a\le
(n-1)/t$, then we can conclude that $\psi K_k\Phi ^*_k\psi _1\in\frak S
_{(n-1)/t,\infty } $.

We have shown in Theorem 3.8 that the norm on operators $K_k$ from a fixed
Sobolev space to another can be estimated in terms of a specific
finite subset
of the symbol seminorms.
Then since the symbol
seminorms of the $K_k$ go to zero rapidly for $k\to\infty $, and the
estimates of the operators $(1-\zeta)\Phi ^*_k\psi _1$ grow at most polynomially in
$k$, the quasinorms ${\bold N}_{(n-1)/t}(\psi K_k\Phi _k^*
\psi _1)$ go to zero rapidly for $k\to\infty $, so the series in 
(4.1) is convergent. This completes the proof. \qed

\enddemo

\proclaim{Corollary 4.2} When $G$ is a $C^\tau $-smooth singular Green
operator of negative order $-t$ and class $0$ on
 a bounded smooth subset 
$\Omega$ of ${\Bbb R}^n$ with boundary  $\Sigma $, then $G\in \frak S_{(n-1)/t,\infty }$.
\endproclaim

\demo{Proof} 
Let $\{\varrho _1,\dots,\varrho
_{J_0}\}$ be a partition of unity subordinate to a relatively open cover $\{U_1,\dots,U_{I_1}\}$ of
$\comega $,
with coordinate maps  $\kappa  _i\colon U_i\to
V_i$ to relatively open
subsets $V_i$ of $ \crnp$, such that any two $\varrho
_{j_1},\varrho _{j_2}$ have support in some $U_i$,
$i=i(j_1,j_2)$. (More on partitions of unity subordinate to covers e.g.\ in \cite{G09},
Lemma 8.4.) We can then write 
$$
G=\sum_{j_1,j_2\le J_0}\varrho _{j_1}G\varrho _{j_2},\tag4.3
$$
where each piece $\varrho _{j_1}G\varrho _{j_2}$ carries over to an
s.g.o.\ $\underline \varrho _{j_1}\underline G\underline \varrho _{j_2}$ on $\rnp$ of the type
 considered in Theorem 4.1, by the coordinate change for $U_i$. 
We apply Theorem 4.1 to the localized pieces. Since
 the coordinate change gives rise to $L_2$-bounded mappings (for
 compact subsets of $U_i$ resp.\ $V_i$ containing $\supp\varrho
 _{j_1}\cup \supp \varrho _{j_2}$ resp.\ $\supp\underline \varrho
 _{j_1}\cup \supp \underline\varrho _{j_2}$), we see
 in view of (2.4) that  $\varrho _{j_1}G\varrho _{j_2}\in \frak
 S_{p,\infty }$ if and only if $\underline \varrho _{j_1}\underline
 G\underline \varrho _{j_2}\in \frak S_{p,\infty }$, with equivalent
 quasinorms ${\bold N}_p$. It follows that the operators
 $\varrho _{j_1}G\varrho _{j_2}$ are in $\frak S_{(n-1)/t,\infty }$, then so
 is their sum by
 (2.2), and the quasinorm can be estimated in terms of the
 quasinorms of the localized pieces.\qed
\enddemo

\subhead 4.2 Asymptotic estimates in selfadjoint cases \endsubhead

Next, we shall show Weyl-type spectral estimates in polyhomogeneous selfadjoint cases.

\proclaim{Theorem 4.3} Let $G$ be a $C^\tau $-smooth polyhomogeneous singular Green operator on
$\rnp$ of order $-t<0$ and class $0$, selfadjoint and $\ge 0$, and let
$\psi (x) =\psi _0(x')\psi
_n(x_n)$, where $\psi _0\in C_0^\infty ({\Bbb R}^{n-1},{\Bbb R})$, 
and $\psi _n\in C_0^\infty ({\Bbb R},{\Bbb R})$
equals $1$ on a neighborhood of $0$. 
Then the eigenvalues of $\psi G\psi $ satisfy
the asymptotic estimate 
$$
\mu _j(\psi G\psi )j^{t/(n-1)}\to c(\psi _0^2 g^0)^{t/(n-1)} \text{ for }j\to\infty ,\tag4.4
$$
where 
$$
 c(\psi _0^2g^0)=\tfrac1{(n-1) (2\pi )^{(n-1)}}\int_{{\Bbb R}^{n-1}}\int_{|\xi '|=1}\tr ((\psi _0(x')^2g^0(x',\xi ',D_n))^{(n-1)/t})\, d\omega dx'.\tag4.5
$$
\endproclaim

\demo{Proof} In addition to the decomposition of $G$ as in Lemma 3.10
using Laguerre operators to the right,
we can decompose by Laguerre operators to the left, writing $G$ in the form 
$$
G=\sum_{l,m\in{\Bbb N}_0}\Phi _lC_{lm}\Phi _m^*, \text{ by defining }C_{lm}=\Phi _l^*G\Phi _m=\Phi _lK_m,\tag4.6
$$
for each $l,m$. For with this definition, 
$$
\sum_{l,m\in{\Bbb N}_0}\Phi _lC_{lm}\Phi _m^*=\sum_{l,m\in{\Bbb
N}_0}\Phi _l\Phi ^*_lG\Phi _m\Phi _m^*=
G,
$$
in view of (3.37). Note that
$$
C_{lm}^*=\Phi _mG^*\Phi _l^*=C_{ml},
$$
since $G=G^*$. 

Each $C_{lm}$ is an operator of  pseudodifferential type on ${\Bbb
R}^{n-1}$, of the form $\Phi _lK_m$ with $K_m$ as defined in Lemma 3.10;
here since we compose with $\Phi _l$ to the left, we get from the
composition rules that $$
\aligned
C_{lm}&=C^0_{lm}+C'_{lm},\quad
C^0_{lm}=\operatorname{OP}'(c_{lm}(x',\xi '))\text{ with }\\
c_{lm}&=\varphi _l(x_n,\sigma
)\circ_n k_m(x',\xi ',D_n)\in C^\tau S^{-t}({\Bbb R}^{n-1}\times {\Bbb R}^{n-1}),\\
C_{lm}'&\colon H^{s-t-\theta }({\Bbb R}^{n-1})\to
H^s({\Bbb R}^{n-1})\text{ for }|s|<\tau , -\tau +\theta <s<\tau.
\endaligned$$
The symbol $c_{lm}(x',\xi ')$ is a polyhomogeneous $C^\tau $-smooth
$\psi $do symbol on ${\Bbb R}^{n-1}$ of order $-t$.

We shall  also consider finite sums
$$
G_M=\sum_{l,m<M}\Phi _lC_{lm}\Phi _m^*.
$$
Defining 
$$
\Cal K_M = \pmatrix \Phi _0 & \Phi _1&\cdots&\Phi
_{M-1}\endpmatrix,\quad
\Cal C_{M}=\bigl(C_{lm}\bigr)_{l,m=0,\dots,M-1}, \tag4.7
$$
a row matrix of Poisson operators, resp.\ an $M\times M$-matrix of
$\psi $do-type operators, we can write
$$
G_M=\Cal K_M \Cal C_M \Cal K_M^*,
$$
it ranges in $\bigoplus_{j<M}H_j$ and vanishes on
$\bigl(\bigoplus_{j<M}H_j\bigr)^\perp$, cf.\ (3.37)-(3.38). Since
$$
G_M=\Cal K_M \Cal K_M^*G \Cal K_M \Cal K_M^*,
$$
 the nonnegativity of $G$ implies that of $G_M$; $\Cal C_M$ is
likewise $\ge 0$.

Denote 
$$
G-G_M=G^\dagger _M,\quad
G^\dagger_M=\sum _{l\text{ or }m\ge M}\Phi _lC_{lm}\Phi _m^*.
$$

The first thing we show is that the contribution from $\psi G^\dagger_M\psi $ to the spectral
asymptotics of $\psi G\psi $ comes from a contributing term where the
constant goes to 0 for
$M\to\infty $, plus a term in a better weak Schatten class. 
Write $\sum_{l\text{ or }m\ge M}=\sum_{l\in {\Bbb N}_0,m\ge
M}+\sum_{l\ge M,m< M}$, and note that
$$
\sum _{l\in{\Bbb
N}_0, m\ge M}\Phi _lC_{lm}\Phi _m^* =\sum _{m\ge M}K_m\Phi _m^*,
$$
and, since $C_{lm}^*=C_{ml}$, 
$$
( \sum _{l\ge M, m< M}\Phi
_lC_{lm}\Phi _m^*)^*=\sum _{l\ge M, m< M}\Phi
_mC_{lm}^*\Phi _l^*=\sum _{m<M, l\ge M}\Phi
_m\Phi _m^*K_l\Phi _l^*=\Psi _M\sum _{l\ge M}K_l\Phi _l^*,
$$
with $\Psi _M=\sum_{m<M}\Phi _m\Phi _m^*=\Cal K_M\Cal K_M^*$ a
selfadjoint projection singular Green operator in $L_2(\rnp)$. Thus
$$
G_M^\dagger=\sum _{l\text{ or }m\ge M}\Phi _lC_{lm}\Phi _m^*=\sum _{m\ge M}K_m\Phi
_m^*+
(\Psi _M \sum_{m\ge M}K_m\Phi _m^*)^*.
$$
 By the calculations in the proof of Theorem 4.1, 
$$
{\bold N}_{(n-1)/t}(\psi\sum _{m\ge M}K_m\Phi
_m^*\psi )\le C \sum_{m\ge M}m^\gamma {\bold N}_{(n-1)/t}(\psi K_m\Phi _m^*
\psi )\to 0 \text{ for }M\to\infty .
$$
This takes care of the first term in $\psi G^\dagger_M\psi $, and for the last
term we insert a cutoff function $\eta $, equal to 1 on a ball containing
$\supp \psi $, leading to:
$$
\psi \Psi _M K_m\Phi _m^*\psi  = \psi \Psi _M \eta K_m\Phi _m^*\varphi  + \eta \psi \Psi _M (1-\eta )K_m\Phi _m^*\psi . 
$$
Here the first term defines a sequence that is rapidly decreasing in $m$ with
respect to the ${\bold N}_{(n-1)/t}$-quasinorm, since $\psi \Psi _M$
has bounded $L_2$-norm independent of $M$.  So the ${\bold
N}_{(n-1)/t}$-quasinorm of the sum over $m\ge M$ goes to 0 for $M\to\infty $. In the second term, $\psi
\Psi _M(1-\eta )$ is of order $-\infty $ and the norm estimates from
$H^{-a}$ to $L_2$ are at most polynomially growing in $M$ for each
fixed $a$, so the summation  in $m$ is in  $\frak S_{p,\infty }$ for
all $p>0$. Thus 
$$
\psi G_M^\dagger \psi = G_{M,1,\psi }^\dagger + G_{M,2,\psi }^\dagger,  \tag4.8
$$
where ${\bold N}_{(n-1)/t}(G_{M,1,\psi }^\dagger)\to 0$ for $M\to \infty $,
and $G_{M,2,\psi }^\dagger\in \bigcap_{p>0}\frak S_{p,\infty }$ for
all $M$. We  view the operators as perturbations of $\psi G_M\psi $, where
Lemma 2.4 $1^\circ$ can be used for $G_{M,2,\psi }$, and Lemma 2.4
$2^\circ$ will later be used for $G_{M,1,\psi }$.

Now consider $\psi G_M\psi $. Let $B'_R$ be a ball containing $\supp
\psi _0$, and let $\zeta\in C_0^\infty ({\Bbb
R}^{n-1})$ be such that $\zeta(x')=1$ on $\overline {B'_R}$. Then
since $\psi _0=\zeta \psi _0$, 
$$
\psi G_M\psi =\sum_{l,m<M} \psi \Phi _l\zeta C_{lm}\Phi _m^*\psi + \sum_{l,m<M} \psi \Phi _l(1-\zeta) C_{lm}\Phi _m^*\psi. 
$$
Here each $\psi \Phi _l(1-\zeta)$ is of order $-\infty $, so $\psi
\Phi _l(1-\zeta) C_{lm}\Phi _m^*\psi$ is continuous from
$H^{-t}(\rnp)$ to $H^a(\rnp)$ for any $a\in{\Bbb R}$, hence is in
$\frak S_{p,\infty }$ for all $p>0$ (note the compactly supported factor $\psi $
to the right). Similarly, we can insert $\zeta$ between $C_{lm}$ and
$\Phi _m^*$, making an error in $\bigcap_{p>0}\frak S_{p,\infty }$, and this leaves
us with having to show the asymptotic estimate for 
$$
\sum_{l,m<M} \psi \Phi _l\zeta C_{lm}\zeta\Phi _m^*\psi.
$$
As a next step, we observe that
$$
(1-\psi _n(x_n))\psi'(x') \Phi _l\zeta C_{lm}\zeta\Phi _m^*\psi
$$
has the factor $(1-\psi _n)\Phi _l$ which is also of order $-\infty $,
so that the operator it enters in,  contributes with a term in
$\bigcap_{p>0}\frak S_{p,\infty }$. Then $\psi _n$ can be omitted from the
operator we have to study.  The same holds to the right. In both cases
we use that a factor $\zeta$ inside the expression assures that the
composed map passes via a compact set.

Now the problem is reduced to studying a sum of terms
$$
\psi _0 \Phi _l\zeta C_{lm }\zeta \Phi _m^*\psi _0.
$$
Here we note that
$$
\psi _0 \Phi _l\zeta C_{lm }\zeta \Phi _m^*\psi _0=[\psi _0 ,\Phi _l]\zeta C_{lm }\zeta \Phi _m^*\psi _0+\Phi _l\psi _0 C_{lm }\zeta \Phi _m^*\psi _0
$$
since $\psi _0\zeta=\psi _0$; here the commutator $[\psi _0,\Phi _l]$
is of 1 step lower order than $\Phi _l$, so the operator is in $\frak
S_{(n-1)/(t+1),\infty }$ and will not contribute to the asymptotic
formula. Doing a similar commutation to the right, we end up with
having to study
$$
G_{M,\psi _0}= \sum_{l,m<M}\Phi _l\psi _0 C_{lm }\psi _0 \Phi _m^*.
$$

With the notation (4.7),
$$
G_{M,\psi _0}= \Cal K_M \Cal C_{M,\psi _0} \Cal K_M^*, \text{ where }\Cal C_{M,\psi _0}=\bigl(\psi _0C_{lm}\psi _0\bigr)_{l,m=0,\dots,M-1}.
$$
The operator is selfadjoint $\ge 0$, so we can use that eigenvalues are
preserved under commutation, to calculate:
$$
\multline
s_j(G_{M,\psi _0})=\mu _j (G_{M,\psi _0})=\mu _j(\Cal K_M \Cal C_{M,\psi
'} \Cal K_M^*)=\\
\mu
_j  ( \Cal C_{M,\psi _0} \Cal K_M^*\Cal K_M)=\mu _j (\Cal C_{M,\psi
'})=s_j (\Cal C_{M,\psi _0}),\endmultline \tag4.9
$$
since $\Cal K_M^*\Cal K_M$ equals the identity matrix in view of
(3.37).

We have now arrived at the consideration of  the operator of  $\psi $do-type
$\Cal C_{M,\psi _0}$. It is of the form 
$$
\Cal C_{M,\psi _0}=\operatorname{OP}'(c_{M,\psi _0})+R,\quad  c_{M,\psi _0}=( \psi _0(x')^2 c_{lm}(x',\xi
'))_{l,m<M},
$$
 where $R\colon H^{-t-\theta }({\Bbb R}^{n-1})\to H^0(B'_{R})$, hence
is in $\frak S_{(n-1)/(t+\theta ),\infty }$ (for some $\theta >0$). A spectral asymptotics
formula is obtained for $\operatorname{OP}'(c_{M,\psi _0})$ by
application of Theorem 2.5, and it extends to $\Cal C_{M,\psi _0}$ by
Lemma 2.4 $1^\circ$, giving:
$$
s_j(\Cal C_{M,\psi _0})j^{t/(n-1)}\to c(c^0_{M,\psi _0})^{t/(n-1)} \text{ for }j\to\infty ,\tag4.10
$$
where the constant is defined as in (2.20) and $c^0_{M,\psi _0}$ is the
principal part of $c_{M,\psi _0}$. Hence $s_j(G_{M,\psi _0})$ likewise
has this behavior (cf.\ (4.9)), and so has $s_j(\psi G_M\psi )$, since
$G_{M,\psi _0}-\psi G_M\psi $ was shown above to be in $\frak
S_{(n-1)/(t+\varepsilon ),\infty }$ for some $\varepsilon >0$.

Since $\tr c^0_{M,\psi _0}(x',\xi ')=\sum_{m<M}\psi _0^2
c^0_{mm}(x',\xi ')=\tr \psi _0^2g^0(x',\xi ',D_n)$ (since these
symbols are related in a similar way as the operators), and a similar
rule holds for traces of powers of these nonnegative selfadjoint operators,
we have that
$$
c(c^0_{M,\psi _0})=c(\psi _0^2g^0_{M})
$$
defined as in (4.5) from the symbol of $\psi _0G_M\psi _0$. So we can now write
$$
s_j(\psi G_M\psi )j^{t/(n-1)}\to c(\psi _0^2g^0_M)^{t/(n-1)}\text{ for }j\to\infty .\tag4.11
$$

Finally we combine this with the information we have on
$G_M^\dagger$. It gives that 
$$
\psi G\psi =\psi G_M\psi +G_{M,1,\psi _0}^\dagger+G_{M,2,\psi _0}^\dagger,
$$
where ${\bold N}_{(n-1)/t}(G_{M,1,\psi }^\dagger)\to 0$ for $M\to \infty $,
and $G_{M,2,\psi }^\dagger\in \bigcap_{p>0}\frak S_{p,\infty }$ for
all $M$. By Lemma 2.4 $1^\circ$, $\psi G_M\psi +G_{M,2,\psi
'}^\dagger$ likewise has the spectral asymptotic behavior
(4.11). Moreover, $c(\psi _0^2g^0_{M})\to c(\psi _0^2g^0)$ for
$M\to\infty $ (using that $ g^0_M\to  g^0$ in the appropriate symbol norms).
Then since $\sup _js_j(G_{M,1,\psi _0}^\dagger)j^{t/(n-1)}\to 0$ for 
$M\to
\infty $, an application of Lemma 2.4 $2^\circ$ gives (4.4) with (4.5). \qed

\enddemo

We have a similar result  in the curved domain situation:

\proclaim{Theorem 4.4} Let $G$ be a selfadjoint nonnegative $C^\tau
$-smooth polyhomogeneous singular Green
operator of negative order $-t$ and class $0$ on
 a bounded smooth subset 
$\Omega$ of ${\Bbb R}^n$ with boundary  $\Sigma $. Then the positive
eigenvalues of $G$ satify an asymptotic estimate$$
\mu _j( G )j^{t/(n-1)}\to c( g^0)^{t/(n-1)} \text{ for }j\to\infty ,\tag4.12
$$
where 
$$
 c(g^0)=\tfrac1{(n-1) (2\pi )^{(n-1)}}\int_{\Sigma }\int_{|\xi '|=1}\tr (g^0(x',\xi ',D_n)^{(n-1)/t})\, d\omega dx'.\tag4.13
$$ 

\endproclaim

\demo{Proof} On a tubular neighborhood of $\Sigma $ we can construct a
family of $C^\infty $-smooth Poisson operators $\widetilde \Phi _l$,
$l\in{\Bbb N}_0$, 
such that their symbol-kernels at $\Sigma $ in  local coordinates have principal
part  equal to $\varphi (x_n,\sigma (\xi '))$ at each $x'\in \Sigma $.

Define 
$$
\aligned
\widetilde C_{lm}&=\widetilde \Phi _l^* G \widetilde \Phi _m\text{ for
}l,m\in{\Bbb N}_0,\\
\widetilde G_M&=\sum_{l,m<M}\widetilde \Phi _l\widetilde C_{lm}\widetilde \Phi
_m^*,\\
\widetilde G_M^\dagger&=G-\widetilde G_M.
\endaligned\tag4.14
$$
Note that $\widetilde C_{lm}^*=\widetilde C_{ml}$, since $G$ is selfadjoint.
We shall show that these operators act in a similar way as in the case
studied in Theorem 4.3.

First consider $\widetilde G_M$. Defining $$
\widetilde{\Cal K}_M = \pmatrix \widetilde \Phi _0 & \widetilde \Phi
_1&\cdots&\widetilde \Phi
_{M-1}\endpmatrix,\quad
\widetilde{\Cal C}_{M}=\bigl(\widetilde C_{lm}\bigr)_{l,m=0,\dots,M-1}, \tag4.15
$$
a row matrix of Poisson operators, resp.\ an $M\times M$-matrix of
$\psi $do-type operators, we can write
$$
\widetilde G_M=\widetilde{\Cal K}_M \widetilde{\Cal C}_M \widetilde{\Cal K}_M^*,
$$
Then since it is selfadjoint $\ge 0$:
$$
s_j(\widetilde G_M)=\mu _j(G_M)=\mu _j(\widetilde{\Cal K}_M \widetilde{\Cal C}_M
\widetilde{\Cal K}_M^*)=\mu _j(\widetilde{\Cal C}_M \widetilde{\Cal
K}_M^*\widetilde{\Cal K}_M ).
$$
Here $\widetilde{\Cal C}_M$ is selfadjoint nonnegative, since for
$\underline u=(u_0,\dots, u_{M-1})\in L_2(\Sigma )^M$,
$$
(\widetilde{\Cal C}_M \underline u,\underline u)_{L_2(\Sigma )^M}=
\sum _{l,m<M}(\widetilde C_{lm}u_m, u_l)_{L_2(\Sigma )}=
\sum _{l,m<M}(G\widetilde \Phi ^*_m u_m,\widetilde\Phi _l^*
u_l)_{L_2(\Omega  )}\ge 0.
$$
Moreover, $\widetilde{\Cal K}_M^*\widetilde{\Cal K}_M$ a selfadjoint nonnegative
$M\times M$-matrix of $C^\infty $-smooth $\psi $do's with principal
part $I_M$ (the identity matrix). Since it is elliptic nonnegative,
its square root $P_M$ is well-defined
(Seeley \cite{S67}). Then we can continue the spectral
calculation as follows: 
$$
\mu _j(\widetilde{\Cal C}_M \widetilde{\Cal
K}_M^*\widetilde{\Cal K}_M )=\mu _j(\widetilde{\Cal C}_MP_M^2)=\mu _j(P_M\widetilde{\Cal C}_MP_M)=\mu _j(\widetilde{\Cal C}_M+R)=s_j(\widetilde{\Cal C}_M+R),
$$
where $R$ is selfadjoint of the form $R=\widetilde{\Cal
C}_M(P_M-I_M)+(P_M-I_M)\widetilde{\Cal C}_MP_M$. Since $\widetilde{\Cal
C}_M\in \frak S_{(n-1)/t,\infty }$, $P_M-I_M\in \frak S_{(n-1),\infty
}$ and $P_M$ is bounded, $R\in \frak S_{(n-1)/(t+1 ),\infty }$ by
(2.3), (2.4).
To $\widetilde{\Cal C}_M$ we 
can apply Theorem 2.4 getting a spectral asymptotics formula, and to
$\widetilde{\Cal C}_M+R$ we can apply the perturbation rule
Lemma 2.4 $1^\circ$ for $s$-numbers. This gives, in view of the
above identifications of eigenvalues:
$$
s_j(\widetilde G_M)j^{t/(n-1)}\to c(c^0_{M})^{t/(n-1)} \text{ for }j\to\infty ,\tag4.16
$$
where the constant is defined as in (2.20) and $c^0_{M}$ denotes the
principal symbol of $\widetilde{\Cal C}_M$.

Now we turn to the study of $\widetilde G_M^\dagger$. Here we use a localization
as in the proof of Corollary 4.2, leading to the study of finitely
many localized pieces
$$
\widetilde G^\dagger_{M,i,j}=\underline\varrho _i(\underline
G-\sum_{l,m<M}\underline{\widetilde \Phi }_l\underline {\widetilde{C}}_{lm}\underline{\widetilde \Phi }_m^*)\underline\varrho _j,
$$
acting in $\rnp$. Consider one of these pieces. We can assume that
$\underline G$ is defined on all of $\rnp$ although it only enters
with cut-off functions around it. Write
$$
\underline G=\sum_{l,m\in{\Bbb N}_0}\Phi _l\underline C_{lm}\Phi ^*_m
,\quad \underline G_M^\dagger=\sum_{l\text{ or }m\ge M}\Phi
_l\underline C_{lm}\Phi ^*_m , 
$$
 as in the proof of Theorem 4.3, then 
$$
\aligned
{\widetilde G}^\dagger_{M,i,j}
&=\underline\varrho _i(\sum_{l,m\in{\Bbb N}_0}\Phi _l\underline
C_{lm}\Phi ^*_m -\sum_{l,m<M}\underline{\widetilde \Phi }_l\underline
{\widetilde{C}}_{lm}\underline{\widetilde \Phi
}_m^*)\underline\varrho _j\\
&=\underline\varrho _i(\underline G_M^\dagger+\sum_{l,m<M}(\Phi
_l\underline C_{lm}\Phi ^*_m-\underline{\widetilde \Phi }_l\underline
{\widetilde{C}}_{lm}\underline{\widetilde \Phi
}_m^*))\underline\varrho _j.
\endaligned
$$
Here
$$
\aligned
\underline\varrho _i(\Phi
_l\underline C_{lm}\Phi ^*_m-\underline{\widetilde \Phi }_l\underline
{\widetilde{C}}_{lm}\underline{\widetilde \Phi
}_m^*)\underline\varrho _j&=\underline \varrho _i((\Phi
_l-\underline{\widetilde \Phi }_l)\underline C_{lm}\Phi
^*_m+\underline{\widetilde \Phi }_l(\underline C_{lm}-\underline
{\widetilde{C}}_{lm}) \Phi
_m^*\\
&\qquad \underline{\widetilde \Phi }_l\underline
{\widetilde{C}}_{lm}(\Phi _m^*-\underline{\widetilde \Phi
}_m^*))\underline\varrho _j,
\endaligned
$$
which is in $\frak S_{(n-1)/(t+\theta ),\infty }$ with $\theta
>0$. This is seen by inspecting each term, e.g.,
$$
\underline \varrho _i(\Phi
_l-\underline{\widetilde \Phi }_l)\underline C_{lm}\Phi
^*_m\underline\varrho _j=
\underline \varrho _i(\Phi
_l-\underline{\widetilde \Phi }_l)\zeta \underline C_{lm}\Phi
^*_m\underline\varrho _j
+\underline \varrho _i(\Phi
_l-\underline{\widetilde \Phi }_l)(1-\zeta )\underline C_{lm}\Phi
^*_m\underline\varrho _j,
$$ with $\zeta \in C_0^\infty ({\Bbb R}^{n-1})$ equal to 1 on $B'_R$, $B_R\supset\supp
 \underline \varrho _i$,
 where the first term is composed of operators in $\frak S_{n-1,\infty
 }$ and $\frak S_{(n-1)/t,\infty}$, and the second term goes from
 $H^{-t}(\rnp)$ to $H^a(B_{R,+})$, any $a>0$. We also use that $\underline C_{lm}-\underline
{\widetilde{C}}_{lm}$ is of order $-t-1$, since the two operators
 have the same principal symbol.

Finally, $\underline\varrho _i\underline G_M^\dagger\underline \varrho
_j$ is of the type already analyzed in the proof of Theorem 4.3;
$$
\underline\varrho _i\underline G_M^\dagger\underline \varrho
_j=\underline\varrho _i\underline G_{M,1}^\dagger\underline \varrho
_j+\underline\varrho _i\underline G_{M,2}^\dagger\underline \varrho
_j,
$$
where $\underline\varrho _i\underline G_{M,2}^\dagger\underline \varrho
_j\in\frak S_{(n-1)/(t+\theta ),\infty }$ and  ${\bold
N}_{(n-1)/t}(\underline\varrho _i\underline G_{M,1}^\dagger\underline \varrho
_j)$ goes to 0 for $M\to\infty $.
Carrying the operators back to $\Omega $ and summing over $i,j$, we find that 
$$
\widetilde G_M^\dagger=\widetilde G_{M,1}^\dagger +\widetilde G_{M,2}^\dagger 
$$
with similar properties. The proof can now be completed in the same
way as the proof of Theorem 4.3.
\qed
\enddemo

For a nonselfadjoint singular Green operator $G$ one gets the result
of Theorem 4.3 when $G^*G$ can be brought on the form $G_1+R$ with
$G_1$ as in Theorem 4.3 and $R$ of order $< -tn/(n-1)$. This holds for
$\tau $ large, in particular when $G$ is $C^\infty $-smooth, and
hereby provides another proof of Theorem 4.10 in \cite{G84}.

\head 5. Results for boundary terms arising from resolvents\endhead

\subhead 5.1 The Dirichlet resolvent \endsubhead

As we saw in Section 4, a severe difficulty in the discussion of spectral
asymptotics formulas for  $C^\tau $-smooth singular Green operators of
order $-t$ is
that although the remainders in compositions are of lower order, say
$-t-\theta $, they only gain an $\frak S_{n/(t+\theta ),\infty
}$-estimate, which is generally not better than the estimate in $\frak S_{(n-1)/t ,\infty
}$ one is aiming for.

In the treatment of singular Green terms coming from resolvents of
differential boundary problems, we shall use another strategy that is based on the special product form
of the s.g.o.-terms. This will even allow nonselfadjointness.

Consider a second-order
strongly elliptic
differential operator $A$ in divergence form
$$
A=-\sum_{j,k=1}^n\partial_ja_{jk}(x)\partial_k +\sum_{j=1}^na_j(x)\partial_j+a_0(x), \tag 5.1
$$
where the $a_{jk}$ and $a_j$ are functions in $H^1_q(\Omega )=\{u\in L_q(\Omega
)\mid Du\in L_q(\Omega )^n\}$ and $a_0\in L_q(\Omega )$, for
some $q>n$. Let $\tau \in \,]0,1-n/q]$; recall that $H^1_q(\Omega
)\hookrightarrow C^\tau (\comega)$. 

In the present section we consider a
smooth $\Omega $; this will be relaxed in Section 6 where there are
additional bounds on $\tau $ adapted to the boundary regularity as in \cite{AGW12}.

It follows from \cite{AGW12} that the solution operator for the nonhomogeneous Dirichlet
problem for $A$ (which exists uniquely, by a variational construction, when a sufficiently large constant has
been added to $A$) maps as follows:
$$
\Cal A=\pmatrix A\\ \gamma _0\endpmatrix \text{ has the inverse }
\Cal B=\pmatrix
R_\gamma & K_\gamma \endpmatrix \colon \matrix
H^{s}(\Omega )\\ \times \\ { H}^{s+\frac32}(\Sigma )\endmatrix\to H^{s+2}(\Omega ),\text{ for } - \tau <s\le 0,\tag5.2
$$
belonging to the nonsmooth $\psi $dbo calculus in the sense that
 $R_\gamma $ and
$K_\gamma $ have the form of  $C^\tau
$-smooth $\psi $dbo's
composed with order-reducing operators, plus  lower-order
remainders. Moreover, $K_\gamma $ extends to a bounded mapping from $H^{s'-\frac12}(\Sigma )$ to
$H^{s'}(\Omega )$ for $s'\in [0,2]$, belonging to the calculus in a
suitable sense.

Furthermore, when $A$ is extended to a uniformly strongly elliptic operator
on ${\Bbb R}^n$ with positive lower bound, the inverse $Q$ there is the sum of a $C^\tau
$-smooth $\psi $do composed with an order-reducing operator, plus
a remainder of lower order. Let $Q_+$ be the truncation to $\Omega $,
i.e., $Q_+=r_\Omega Qe_\Omega $, where $r_\Omega $ restricts to
$\Omega $ and $e_\Omega $ extends by zero. Then the resolvent can be written
$$
R_\gamma =Q_++G_\gamma ,\quad G_\gamma =-K_\gamma \gamma _0Q_+.\tag5.3
$$

Here the product structure of $G_\gamma
$, passing via the $(n-1)$-dimensional manifold $\Sigma $, will be useful in Schatten class considerations. 

The choice of the space $H^1_q$ for the coefficients in (5.1) is
governed by the fact that a result of Marschall in \cite{M88} shows
that this allows a comparison with the $x$-form operator
$A^\times=-\sum_{j,k=1}^na_{jk}(x)\partial_j\partial_k $,
assuring that $A-A^\times$ is continuous from $H^{s+2-\theta }(\Omega )$
to $H^s(\Omega )$ for $-\tau +\theta <s\le 0$.

Our strategy will be to approximate $G_\gamma $ by a composition of
$C^\infty $-coefficient operators, obtained by approximation the original boundary value problem
by  $C^\infty $-coefficient problems, where the desired estimates are
well-known.

\proclaim{Proposition 5.1} There is a sequence of $C^\infty $-smooth invertible $\psi
$dbo systems $\Cal A_k=\pmatrix A_k\\\gamma
_0\endpmatrix$, with inverses $A_k^{-1}=Q_k$ on ${\Bbb R}^n$,
$\Cal A_k^{-1}=\Cal B_k=\pmatrix R_{\gamma ,k}& K_{\gamma
,k}\endpmatrix$ on $\Omega $, such that
the symbol seminorms of $A-A_k$ converge to $0$ in 
$H^1_qS^2({\Bbb R}^n\times {\Bbb R}^n)$, and in 
$C^{\tau
'}S^2({\Bbb R}^n\times {\Bbb R}^n)$ for all $\tau '\in\,]0,\tau [\,$,
and when $k\to\infty $,$$
\aligned
\|A -A_{k}\|_{\Cal L(H^{s+2}({\Bbb R}^n ),
H^{s}({\Bbb R}^n ))}&\to 0\text{ for each }|s|<\tau ,\\
\|\Cal A-\Cal A_k\|_{\Cal L(H ^{s+2}(\Omega ), H^s(\Omega )\times
H^{s+\frac32}(\Sigma ))}&\to 0
\text{ for each } |s|<\tau ,\\
\|Q -Q_{k}\|_{\Cal L(H^s({\Bbb R}^n ),
H^{s+2}({\Bbb R}^n ))}&\to 0\text{ for each }s\in \,]-\tau ,0] ,\\
\|\Cal B -\Cal B_{k}\|_{\Cal L(H^s(\Omega )\times H^{s+\frac32}(\Sigma ),
H^{s+2}(\Omega ))}&\to 0\text{ for each }s\in \,]-\tau ,0] ,\\
\|K_\gamma -K_{\gamma ,k}\|_{\Cal L( H^{s'-\frac12}(\Sigma ),
H^{s'}(\Omega  ))}&\to 0\text{ for each } s'\in [0,2] .
\endaligned\tag5.4
$$
 
\endproclaim

\demo{Proof}
As accounted for in Section 2,  $A_k=\varrho _k*A$ 
converges to $A$ in $H^1_qS^2$ and in  $C^{\tau '}S^2$ for any $\tau '\in \,]0,\tau
[\,$. Starting with a sufficiently large $k$ and relabelling, 
we can assure that the operators $A_k$ are uniformly strongly elliptic
on ${\Bbb R}^n$, uniformly in $k$. (Note that the $A_k$ are differential operators and
only their coefficients  are modified by convolution with $\varrho _k$.)
We can assume that a sufficiently large constant has
been added to $A$ such that all the operators have a lower
bound greater than a positive constant. 

By Theorem 2.3 ff., the first two lines in (5.4) are valid. For the
second line we note that for
differential operators, the restriction to functions on $\Omega $ is
tacitly understood, and that $\gamma _0$
cancels out.

 The other statements will be obtained
by following each step of the construction in
\cite{AGW12}. First there is the
construction in Sect.\ 4.2 there, where parametrices are constructed from
the principal symbols in a localized situation. The principal interior
symbol and boundary symbol
are determined pointwise in $x$ or $x'$ by the standard
theory; when $C^{\tau '}$-dependence is taken into account, one
similarly gets
the $C^{\tau '}$-dependence for the inverted symbols. Then the
convergence of the principal symbol of $A_k$ to that of $A$ in
$C^{\tau '}$ implies a convergence in $C^{\tau '}$ of the inverse
principal symbol and  boundary symbol. By Theorem 2.3 ff. and Theorem
3.8, the resulting operator families converge in the asserted Sobolev norms;
note also that the Sobolev operator norms are bounded 
uniformly in $k$. 

Next, we go to the $\lambda $-dependent construction in Sect.\ 4.3 of \cite{AGW12},
applied to the operators on $\Omega $. Here the remainder in the
calculation of $\Cal A_k(\lambda )\Cal B^0_k(\lambda )-I$ will have
not only lower order, but also a small operator-norm in terms of $\lambda
$, uniformly in $k$.

It is found that for sufficiently large $\lambda $ in a sector around
${\Bbb R}_-$ the exact inverses are determined by Neumann series
arguments and have the form of the sum of a
parametrix coming from the calculus and a lower-order term. 
The convergence in the third and fourth
lines of (5.4) follows from the estimates, for the values of $s$ allowed in the inverse construction for the nonsmooth problem, namely for $s\in \,]-\tau ,0]$.

For the last statement in (5.4) we use the additional information
worked out in Sect.\ 5.2 of \cite{AGW12}, again examined with a view to
the approximation in $C^{\tau '}$ for $k\to\infty $.  
\qed
\enddemo

\proclaim{Theorem 5.2} 

$1^\circ$ The $\psi $do term $Q_+$ in {\rm (5.3)} has the spectral
behavior
$$
s _j( Q_+  )j^{2/n}\to c( q^0_+)^{2/n} \text{ for }j\to\infty ,\tag5.5
$$
where the constant is defined from the principal symbol $q^0$  (or $a^0=(q^0)^{-1}$):
$$
 c(q^0_+)=\tfrac1{n(2\pi )^{n}}\int_{\Omega  }\int_{|\xi
|=1} (q^{0*}(x,\xi )q^{0}(x,\xi ))^{n/4}\,d\omega dx
=\tfrac1{n(2\pi )^{n}}\int_{\Omega  }\int_{|\xi
|=1} (a^{0*}a^{0})^{-n/4}\,
d\omega dx.\tag5.6 
$$ 

$2^\circ$ The singular Green term $G_\gamma
=-K_\gamma \gamma _0Q_{+}$ (cf.\ {\rm (5.3)}) in the solution operator
to the Dirichlet problem for  $A $ has the spectral behavior
$$
s _j( G_\gamma  )j^{2/(n-1)}\to c( g^0)^{2/(n-1)} \text{ for }j\to\infty ,\tag5.7
$$
where the constant is defined from the principal symbol $g^0$ of $G_\gamma $:
$$
 c(g^0)=\tfrac1{(n-1)  (2\pi )^{(n-1)}}\int_{\Sigma }\int_{|\xi
'|=1}\tr ((g^{0*}(x',\xi ',D_n)g^{0}(x',\xi ',D_n))^{(n-1)/4})\,
d\omega dx',\tag5.8 
$$ 
cf.\ also {\rm (5.11)} below.
\endproclaim

\demo{Proof} We approximate $Q$ and $K_\gamma $ by the sequences $Q_k$
and $K_{\gamma ,k}$ as in Proposition 5.1. It is known in the
$C^\infty $-case that the exact
solution operators $Q_k$, $\Cal B_k=\pmatrix R_{\gamma ,k}&K_{\gamma
,k}\endpmatrix$ belong to the $\psi $dbo calculus (cf.\ e.g.\
\cite{G96} or the recent account in \cite{G09}, Sect.\ 11.3). 

$1^\circ$. The statements (5.5)--(5.6) are well-known for $Q_{k,+}$, 
see e.g.\ \cite{G96}
Proposition 4.5.3. It follows from the third line in (5.4) that
$Q_+-Q_{k,+}\to 0$ in the norm of $\Cal L(H^0(\Omega ), H^2(\Omega ))$,
hence in the quasinorm of $\frak S_{n/2,\infty }(L_2(\Omega ))$. Then
the result follows by an application of Lemms 2.4 $2^\circ$.

$2^\circ$. The
statements (5.7)--(5.8) for $G_{\gamma ,k}=-K_{\gamma ,k}\gamma
_0Q_{k,+}$ are known for these operators from \cite{G84}.

Now we have 
$$
G_\gamma -G_{\gamma ,k}=(-K_\gamma +K_{\gamma ,k})\gamma _0Q_+-K_{\gamma ,k}(\gamma _0Q_+-\gamma _0Q_{k,+}).
$$
  Here $\gamma _0Q_+$ is bounded from $H^0(\Omega )$ to   $H^{\frac32}(\Sigma )$. Writing it as $\Lambda _0^{-\frac32}\Lambda  _0^{\frac32}\gamma _0Q_+$, where the $\Lambda _0^r$ are
  homeomorphisms from $H^s(\Sigma )$ to $H^{s-r}(\Sigma )$ for all $s$,
  with $(\Lambda _0^{r})^{-1}=\Lambda _0^{-r}$, we see that $\gamma
  _0Q_+$ belongs to $\frak  S_{(n-1)/(3/2),\infty }(L_2(\Omega ),L_2(\Sigma ))$, as the composition of the
  operator $\Lambda _0^{-\frac32}\in \frak  S_{(n-1)/(3/2),\infty } (L_2(\Sigma ))$
  and the bounded operator $\Lambda  _0^{\frac32}\gamma _0Q_+$ from
  $L_2(\Omega )$ to $L_2(\Sigma )$. 
Similarly, since $K_\gamma  -K_{\gamma ,k}= (K_\gamma  -K_{\gamma ,k})\Lambda
  _0^{\frac12}\Lambda  _0^{-\frac12}$, where $\Lambda _0^{-\frac12}\in
  \frak S_{(n-1)/(1/2),\infty }(L_2(\Sigma ))$, and the norm of $(K_\gamma  -K_{\gamma ,k})\Lambda
  _0^{\frac12}$ from $L_2(\Sigma )$ to
$L_2(\Omega )$ goes to 0 in view of (5.4),   
$$
{\bold N}_{(n-1)/(1/2)}(K_\gamma -K_{\gamma ,k})\to 0 \text{ for }k\to \infty .\tag5.9
$$
It follows  by (2.3) that
$$
{\bold N}_{(n-1)/2}((K_\gamma -K_{\gamma ,k})\gamma _0Q_+)\to 0 \text{
for }k\to \infty . 
$$

For the other term, note that ${\bold N}_{(n-1)/(1/2)}(K_{\gamma ,k})$ is bounded in $k$,
and \linebreak${\bold N}_{(n-1)/(3/2)}(\gamma _0Q_+ -\gamma _0Q_{k,+})\to 0
\text{ for }k\to \infty $ by (5.4), so also
$$
{\bold N}_{(n-1)/2}(K_{\gamma ,k}(\gamma _0Q_+-\gamma _0Q_{k,+}))\to 0 \text{
for }k\to \infty . 
$$

Together these statements show that
$$
{\bold N}_{(n-1)/2}(G_\gamma -G_{\gamma ,k})\to 0 \text{ for }k\to \infty .
$$
We can then apply Lemma 2.4 $2^\circ$ to the decomposition $G_\gamma
=G_{\gamma ,k}+(G_\gamma -G_{\gamma ,k})$, concluding the assertion in
the theorem, since $g^0_k\to g^0$ for $k\to\infty $.
\qed

\enddemo

\example{Remark 5.3}
Because of the special structure of $G_\gamma $ it is easy to
determine the constant, in the following way: When all operators are smooth we have,
denoting $\gamma _0Q_+=T$:
$$
s_j(G_  \gamma )^2=\mu _j(G_\gamma ^* G_\gamma )=\mu _j(T^*K^*_\gamma K_\gamma
T
)=\mu _j( 
K^*_\gamma K_\gamma TT^*)=\mu _j( P_1P_1' )=\mu _j(P_2 (P'_2)^2P_2),
$$ where $P_1=K^*_\gamma K_\gamma $, $P_1'=TT^*$ are selfadjoint
nonnegative elliptic operators with square roots $P_2$ resp.\ $P'_2$. Then
$$
s_j(G_\gamma )=s_j(P_2'P_2).
$$
To calculate the principal symbols at a point $x'\in\Sigma $ we
consider the principal symbol of $A$ in local
coordinates for $|\xi '|\ge 1$: 
$$
\align
a^0(x',\xi ',\xi _n)&=s_0(x')(\xi _n-\lambda ^+(x',\xi '))(\xi
_n-\lambda ^-(x',\xi '))\\
&=s_0(x')(i\xi _n+\kappa  ^+(x',\xi '))(-i\xi
_n+\kappa  ^-(x',\xi ')),
\endalign
$$
where $\operatorname{Im}\lambda ^\pm\gtrless 0$, and $\kappa ^+=-i\lambda
^+$, $\kappa ^-=i\lambda ^-$ have positive real part. Then $Q$ has principal symbol
$$
q^0=\frac1{a^0}=\frac{1}{s_0(\kappa ^++\kappa ^-)}\Bigl(\frac1{\kappa ^++i\xi _n}+\frac1{\kappa ^--i\xi _n}\Bigr),
$$
 and it follows from the rules of calculus for $\psi $dbo's that $T$
 and $K_\gamma $ 
have principal symbols
$$
t^0=\frac1{s_0(\kappa ^++\kappa ^-)(\kappa ^--i\xi _n)},\quad k^0=\frac1{\kappa ^++i\xi _n}.\tag5.10
$$ Then 
$P'_1=TT^*$ and $P_1=K^*_\gamma K_\gamma $ have principal symbols 
$$
p_1^{\prime 0}=\frac1{|s_0(\kappa ^++\kappa
^-)|^{2}2\operatorname{Re}\kappa ^-}, \quad p_1^0=\frac1{2\operatorname{Re}\kappa ^+}, 
$$
and (cf.\ Theorem 2.5)
$$
c(g^0)=c(p_2^{\prime 0}p_2^0)=\tfrac1{(n-1) (2\pi )^{(n-1)}}\int_{\Sigma }\int_{|\xi '|=1}(4|s_0(\kappa ^++\kappa
^-)|^{2} \operatorname{Re}\kappa ^-\operatorname{Re}\kappa ^+ )^{-(n-1)/4}\, d\omega dx'.\tag5.11
$$
The formula carries over to the nonsmooth case since the symbols
converge in $C^{\tau '}$ ($\tau '<\tau $) for $k\to\infty $.
\endexample

\subhead 5.2 Resolvent differences \endsubhead

Next, we consider the difference between the Dirichlet resolvent and the
resolvent of a Neumann problem for $A$ defined by a boundary condition
$\chi u-C\gamma _0u=0$ as studied in \cite{AGW12}, Sect.\ 7. Here
$\chi u$ is the conormal derivative,
$$
\chi u=\sum_{j,k=1}^nn_j\gamma _0a_{jk}\partial_ku,\tag 5.12
$$
where $\vec n=(n_1,\dots,n_n)$ is the interior unit normal on $\Sigma $.
Moreover, $C$ is a first-order tangential differential operator $C=c\cdot
D_\tau +c_0$, where $c=(c_1,\dots,c_n)$, all $c_j\in
H^1_q(\Sigma )$ ($D_\tau =-i\partial_\tau $ as defined in
\cite{AGW12}, Sect.\ 2.4). 
To assure invertibility, let us assume that the sesquilinear form
$$
a_{\chi
,C}(u,v)=\sum_{j,k=1}^n(a_{jk}\partial_ku,\partial_jv)_{L_2(\Omega )} 
+(\sum_{j=1}^na_j\partial_ju+a_0u,v)_{L_2(\Omega )} +(C\gamma _0u,\gamma _0v)_{L_2(\Sigma  )}  \tag 5.13
$$
satisfies $\operatorname{Re}a_{\chi ,C}(u,u)\ge c_0\|u\|^2_{H^1(\Omega
)}-k\|u\|^2_{L_2(\Omega )}$ for $u\in H^1(\Omega )$; then we can add a
constant (absorbed in $a_0$) such that the realization $A_{\chi ,C}$
in $L_2(\Omega )$ defined from $a_{\chi ,C}$ by variational theory has
its spectrum in a sector in $\{\operatorname{Re}z>0\}$.
As shown in \cite{AGW12} there holds a Krein-type resolvent formula
$$
A_{\chi ,C}^{-1}-A_\gamma ^{-1}=K_\gamma L^{-1}(K'_\gamma )^*;\tag5.14
$$
here $L$ acts like $C-P_{\gamma ,\chi }$, where $P_{\gamma ,\chi }$ is
the Dirichlet-to-Neumann operator.
The boundary value problem is elliptic, hence so is $L=C-P_{\gamma
,\chi }$, so $D(L)=H^{\frac32}(\Sigma )$ (cf.\ \cite{AGW12} Th. 7.2), mapping this space
homeomorphically onto $H^{\frac12}(\Sigma )$. 

It is shown in Sect.\ 5 of \cite{AGW12} that $$
P_{\gamma ,\chi }\colon H^{s-\frac12}(\Sigma )\to H^{s-\frac32}(\Sigma )$$
is continuous for $s\in [0,2]$, and since multiplication by the
coefficients $c_j$ 
in $C$ preserves $H^s(\Sigma )$ for $|s|\le 1$ (cf.\ \cite{AGW12}
(2.29)), $C$ has this mapping property for $s\in [\frac12,2]$. So $L$
extends to a continuous map of this kind for $s\in [\frac12,2]$.

When we now also approximate $C$ by smoothed out operators $C_k$, and
$P_{\gamma ,\chi }$ is approximated by smooth $\psi $do's $P_{\gamma
,\chi ,k}$ as a result of the earlier mentioned smoothing of $A$, we
get a sequence $L_k=C_k-P_{\gamma ,\chi ,k}$ such that
$$
\|L-L_k\|_{{\Cal L}(H^{s-\frac12}(\Sigma ), H^{s-\frac32}(\Sigma
))}\to 0\text{ for }k\to\infty , \text{ all }s\in [\tfrac12,2].
$$
Since $L$ is an invertible operator from $H^{\frac32}(\Sigma )$ to
$H^{\frac12}(\Sigma )$, $L_k$ is likewise so for large $k$, say $k\ge
k_0$, by a Neumann series argument. Moreover, the inverses $L_k^{-1}$
are estimated uniformly in $k\ge k_0$, and 
$$
\|L^{-1}-L_k^{-1}\|_{{\Cal L}(H^{\frac12},
H^{\frac32})}=\|L^{-1}(L_k-L)L_k^{-1}\|_{{\Cal
L}(H^{\frac12}, H^{\frac32})}\to 0\text{ for }k\to\infty .
$$
We can then show:

\proclaim{Theorem 5.4} The singular Green term $G_C
=K_\gamma L^{-1}(K_\gamma ')^*$ in the Krein resolvent formula {\rm
(5.14)}
  has the spectral behavior
$$
s _j( G _C)j^{2/(n-1)}\to c( g^0_C)^{2/(n-1)} \text{ for }j\to\infty ,\tag5.15
$$
where the constant is defined from the principal symbol $g^0_C$ of $G _C$
by formula {\rm (5.8)}, cf.\ also {\rm (5.17)} below.

\endproclaim 
    
\demo{Proof}
Define $G_{C,k}=K_{\gamma ,k}L_k^{-1}(K_{\gamma ,k}')^*$, and write
$$
G_C-G_{C,k}=(K_\gamma -K_{\gamma ,k})L^{-1}(K_{\gamma }')^*+K_{\gamma ,k}(L^{-1}-L_k^{-1})(K_{\gamma }')^*+K_{\gamma ,k}L_k^{-1}((K_{\gamma }')^*-(K_{\gamma ,k}')^*).\tag5.16
$$
We shall show that all three terms have $\bold
N_{(n-1)/2}$-quasinorms going to 0 for $k\to\infty $. Then since the
statements in the theorem are well-known for $G_{C,k}$, the result follows
for $G_C$ by application of Lemma 2.4 2$^\circ$.

Since $K'_\gamma $ is of the same kind as $K_\gamma $, its adjoint
satifies for all $s'\in [0,2]$ that
$(K'_\gamma )^*\in \Cal L(H^{-s'}_0(\comega), H^{-s'+\frac12}(\Sigma
))$, and $(K'_\gamma )^*-(K'_{\gamma ,k})^*$ goes to zero in these
operator norms for $k\to\infty $.
 
In the first term in (5.16), the boundedness of $(K'_\gamma )^*$ from
$H^0(\Omega )$ to $H^\frac12(\Sigma )$ and that of  $L^{-1}$ from
$H^\frac12(\Sigma )$ to $H^{\frac32}(\Sigma )$ imply that
$L^{-1}(K'_\gamma )^*\in \Cal L(H^0(\Omega ), H^{\frac32}(\Sigma ))$;
hence it lies in $\frak S_{(n-1)/(3/2),\infty }$ (by an argumentation
as in the proof of Theorem 5.2). Together with
(5.9) this implies that the composed operator goes to 0 in $\frak
S_{(n-1)/2,\infty }$ for $k\to\infty $.

For the second term, we use that the norm of
$(L^{-1}-L_k^{-1})(K'_\gamma )^*$ in $\Cal L(H^0(\Omega
),H^{\frac32}(\Sigma ))$ goes to zero, and for the last term 
we use that the norm of
$L_k^{-1}((K'_\gamma )^*-(K_{\gamma ,k}')^*) $ in $\Cal L(H^0(\Omega
),H^{\frac32}(\Sigma ))$ goes to zero; both are combined with the fact that
$K_\gamma \in \frak S_{(n-1)/(1/2),\infty }$. This completes the proof.\qed
\enddemo

\example{Remark 5.5} Also here, the constant can be determined from a
formula with $\psi $do's on $\Sigma $. In local coordinates, $K_\gamma
$ has the principal symbol (5.11), and since $A'$ has the principal
symbol $\overline{a}^0$, $K'_\gamma $ has the principal symbol
$(\overline {\kappa ^-}+i\xi _n)^{-1}$, 
and that of $(K'_\gamma )^*$
 is $( \kappa ^--i\xi _n)^{-1}$. Let $P_1''=(K'_\gamma )^*K'_\gamma $,
 with square root $P_2''$; here
 $p_2^{\prime\prime0}=(2\operatorname{Re}\kappa ^-)^{-1/2}$. Then, with
 notation from Remark 5.3 (recall $p_2^{0}=(2\operatorname{Re}\kappa ^+)^{-1/2}$), one has in case of smooth operators:
$$
s_j(G_C)^2=\mu _j(K'_\gamma(L^*)^{-1}K_\gamma ^*K_\gamma
L^{-1}(K'_\gamma )^*) 
=\mu _j((L^*)^{-1}P_2^2 L^{-1}(P_2'')^2)= s_j(P_2 L^{-1}P_2'')^2, 
$$
and hence, with $l^0$ denoting the principal symbol of $L$,
$$
 c(g^0_C)=\tfrac1{(n-1) (2\pi )^{(n-1)}}\int_{\Sigma }\int_{|\xi '|=1}(4(l^0)^2\operatorname{Re}\kappa ^+\operatorname{Re}\kappa ^-)^{-(n-1)/4} \, d\omega dx'.\tag5.17
$$ 
The formula extends to nonsmooth cases by the approximation used in 
 Theorem 5.4.

\endexample

\proclaim{Corollary 5.6} The operators $A_\gamma ^{-1}$ and $A_{\chi
,C}^{-1}$ have the same spectral asymptotic behavior as $Q_+$ in
Theorem {\rm 5.2}.
\endproclaim

\demo{Proof} The statement for $A_\gamma ^{-1}$ is found by
application of Lemma 2.4 $1^\circ$ and Theorem 5.2 to (5.3). Next, the
statement for $A_{\chi ,C}$ follows similarly from (5.14) and Theorem 5.4.\qed 
\enddemo

\example{Remark 5.7}
The estimates extend to differences between resolvents at arbitrary points
of $\varrho (A_\gamma )\cap \varrho (A_{\chi ,C})$ in the same way as
in \cite{G11}, using (3.7) there.
\endexample

\head 6. Nonsmooth boundaries \endhead

In \cite{AGW12} also open sets with nonsmooth
boundaries were allowed. While the
coefficients
$a_{jk}, a_j$ were taken in $H^1_q(\Omega )$ ($a_0\in L_q(\Omega )$)
as above, the domain $\Omega $ was taken to be
$B^{\frac32}_{p,2}$ (meaning that the boundary is locally the graph of
a function in $B^{\frac32}_{p,2}({\Bbb R}^{n-1})$), possibly with different choices of
$p$ and $q$.
More precisely it was required that $p,q>2$, $p<\infty $, with
$$
1-\tfrac nq\ge \tfrac12 -\tfrac{n-1}p\equiv \tau _0>0;\tag6.1
$$
then $\tau $ was taken as
$\tau _0$.
This domain regularity includes $C^{\frac32+\varepsilon }$ for small
$\varepsilon $, and is included in $C^{1+\tau _0}$. 

The interest of the Besov space $B^{\frac32}_{p,2}$ is that it hits the
differentiability index $\frac32$ in an exact way, better that
$C^{\frac32+\varepsilon }$. It comes in as the boundary value space
(trace space) for
the anisotropic Sobolev space $$
W^{2}_{(2,p)}(\rnp)=\{f\mid
D^\alpha f\in L_2(\rp;L_p({\Bbb R}^{n-1})), |\alpha |\le 2
\}.\tag6.2$$ 

\proclaim{Proposition 6.1} $1^\circ$ Let $0<\delta <\delta _1<1$. When $\Omega $ is a bounded
$C^{1+\delta _1}$-domain, there exists a sequence of  bounded $C^\infty $-domains
$\Omega _l$ with $C^{1+\delta  }$-diffeomorphisms $\lambda 
_l\colon {\comega  }\simto \comega_l
$
(sending $\partial\Omega =\Sigma $ onto $\partial\Omega _l=\Sigma _l$),
such that $\lambda  _l\to \operatorname{Id}$ in $C^{1+\delta }$ on a
neighborhood $\Omega '$ of $\comega$,
 and in this sense $\comega _l\to \comega $,
for $l\to\infty $.

$2^\circ$ When $\Omega $ is a bounded domain with
$B^{\frac32}_{p,2}$-boundary, it is $C^{1+\tau _0}$-smooth, and there 
exists a sequence of  bounded $C^\infty $-domains
$\Omega _l$ with $C^{1+\tau _0 }$-diffeomorphisms $\lambda 
_l\colon {\comega  }\simto \comega_l$  (associated to the $B^{\frac32}_{p,2}$-boundary structure as in \cite{AGW12})
such that $\lambda  _l\to \operatorname{Id}$ in $C^{1+\tau _0 }$ on a
neighborhood $\Omega '$ of $\comega$,
 and in this sense $\comega _l\to \comega $,
for $l\to\infty $. Moreover, the estimates for
the associated operators
studied in \cite{AGW12} for second-order problems hold boundedly in $l$.

\endproclaim

\demo{Proof} $1^\circ$.
The $C^{1+\delta _1}$-structure is defined by a family
of $C^{1+\delta _1}$-charts $\kappa _j\colon V_j\to U_j$, $j=1,\dots,N$, where $V_j$ and $U_j$ are
open subsets of ${\Bbb R}^n$ and  $\kappa _j$ is a $C^{1+\delta _1}$-diffeomorphism 
from $V_j$ to $U_j$, sending $V_{j,+}=V_j\cap \crnp$ resp.\ $V_{j,0}=V_j\cap ({\Bbb
R}^{n-1}\times \{0\})$ onto $U_{j,+}=U_j\cap \comega$ resp.\ $U_{j,0}=U_j\cap \partial\Omega $, and where $U\equiv \bigcup _{j\le N}U_j\supset \comega$.
We can choose bounded open subsets $V'_j, V''_j$, mapped onto $U'_j,U''_j$, such that
$\overline V''_j\subset
V'_j\subset \overline V'_j\subset V_j$,
$\overline U''_j\subset
U'_j\subset \overline U'_j\subset U_j$, 
with $U''\equiv\bigcup _{j\le N}U''_j\supset \comega$, and there is a smooth partition of unity 
$\{\psi _j\}_{j\le N}$ for $\overline U''$ with $\supp \psi _j\subset U'_j$,
$\psi _j\ge c_j>0$ on $U''_j$ for each $j$.
With $\varrho _l$ denoting an
approximative unit with $\supp \varrho _l\subset  B_{1/l}$,
the convolution $\varrho _l* \kappa _j$ is well-defined on
$V_{j,l}=\{x\in V_j\mid \dist(x,\partial V_j)>1/l\}$, and 
$$
\varrho _l* \kappa _j\to \kappa _j \text{ in }C^{1+\delta
}(
V'_j)\text{ for }l\to\infty ,\; l\ge l_0,
$$
with $l_0$ so large that $\overline V'_j\subset V_{j,l_0}$ for all
$j\le N$.
Then define the mapping $\lambda _l$ on the neighborhood $U=\bigcup
U_j$ of $\comega$ by
$$
\lambda _l=\sum_{j\le N}(\varrho _l*\kappa _j)\circ(\psi _j\kappa _j^{-1});\tag6.3
$$
it converges to the identity on $\overline U''$, in the $C^{1+\delta
}$-norm.

Let $\Omega _l=\lambda _l(\Omega )$; for large $l$ it is contained in
a small neighborhood of $\comega $. The mappings $\varrho _l*\kappa
_j$ from $V'_j$ to neighborhoods of $\overline U'_j$ are $C^\infty $ and are
bijective from the sets $V''_j$ to their images $\widetilde U_{j,l}$
in the neighborhoods of  $\overline U''_j$, for $l$
sufficiently large. Note that they are $C^\infty $ from the subset
$V''_{j,+}$ (with the flat boundary piece $V''_{j,0}\subset {\Bbb R}^{n-1}\times
\{0\}$) to $\widetilde U_{j,l}\cap \comega_l$, so they provide a system of
$C^\infty $-charts
for $\comega_l$. This is then indeed a $C^\infty $-smooth subset of ${\Bbb R}^n$.

Since $\lambda _l\to \operatorname{Id}$ on $\overline U''$, it is a
$C^{1+\delta }$-diffeomorphism from $U''$ to $\lambda _l(U'')$ for $l$
sufficiently large. In particular, it then defines a $C^{1+\delta
}$-diffeomorphism from $\comega$ to $\comega_l$, converging to the
identity for $l \to\infty $; in this sense $\comega
_l\to\comega$ for $l\to\infty $.

$2^\circ$. In the case of a domain with $B^{\frac32}_{p,2}$-boundary, one
can likewise approximate with $C^\infty $-domains. It is shown in
\cite{AGW12} that such domains are $C^{1+\tau _0}$, so
one can use $1^\circ$ to do a
$C^{1+\tau }$-approximation for any $\tau <\tau _0$. One can also do a
more precise construction
adapted to the structure explained in \cite{AGW12}:

We refer to the notation there, see in particular Lemma
2.2, Prop.\ 2.7 and Remark
2.9 there. Recall that $\comega$ is covered by patches $U_0,\dots ,
U_J$, where $\overline U_0\subset \Omega $, and the other $U_j$, together covering $\Sigma $, 
 are
such that $U_j\cap \Omega $ is of the form $U_j\cap \{x\mid x_n>\gamma
_j(x')\}$ for a function $\gamma _j\in B^{\frac32}_{p,2}({\Bbb
R}^{n-1})$ (after a rotation of coordinates). From $\gamma _j$ one
constructs (in Prop.\ 2.7 of \cite{AGW12}) a $C^{1+\tau
_0}$-diffeomorphism $F_j$ on ${\Bbb R}^n$ mapping $\rnp$ to $ \{x\mid x_n>\gamma
_j(x')\}$.
For each of these $U_j$, we can smoothen out the
boundary defining function $\gamma _j(x')$ by convolution with an
approximate identity $\varrho _l(x')$; then $\varrho _l*\gamma _j$ 
converges to $\gamma _j$ in
$B^{\frac32}_{p,2}$-norm for $l\to\infty $. 
The associated diffeomorphisms  $F_{j,l}$ 
constructed as in Prop.\ 2.7  are $C^\infty $ (since the
lifting used in Lemma 2.2 of \cite{AGW12} lifts $C^\infty $ to $C^\infty $), and converge to $F_j$
in $C^{1+\tau _0}$. For $U _0$ we simply use the identity
diffeomorphism for all $l$. As under $1^\circ$, the approximations can
be pieced
together by a suitably chosen partition of unity, leading to an
approximation of $\Omega $ by $C^\infty $-domains $\Omega _l$, such
that the boundary converges in $B^{\frac32}_{p,2}$ and the
diffeomorphisms converge in $C^{1+\tau _0}$. Moreover, the various other
constructions in \cite{AGW12} can be followed through these coordinate
changes.
\qed

\enddemo

For $C^1$-domains it is not hard to prove that $1^\circ$ holds with
approximation in $C^1$.

Let $\Omega _l=\uomega$ be one of the approximating domains, with
a $C^{1+\tau }$-diffeomorphism $\lambda \colon \comega\to \ucomega$. 
When $u$ is a function on $\comega $, denote by $\uu$ the function
on $\overline{\uomega }$ defined by 
$$
\uu(y)=u(\lambda ^{-1}(y))\text{ for }y\in \overline{\uomega }, \text{ also
denoted }
u\circ \lambda ^{-1}\text{ or }(\lambda ^{-1})^*u.\tag6.4
$$
An operator $B$ in $L_2(\Omega )$ carries over to the operator $\uB$
in $L_2(\uomega)$ by
$$
\uB\, \uu = \underline {(Bu)},\text{ also written as }(\lambda ^{-1})^* B \lambda ^*\uu. \tag6.5
$$
We have by the rules
for coordinate changes:
$$
(Bu,v)_{L_2(\Omega )}=\int_{\Omega }Bu(x)\bar v(x)\,dx=\int_{\uomega }\uB\,\uu(y)\bar {\uv}(y)J(y)\,dy=(\uB\,\uu,J\bar{\uv})_{L_2(\uomega)},
$$
 where $J$ denotes the Jacobian $|\det (\partial (\lambda^{-1})
 _j/\partial y_k)_{j,k=1,\dots ,n}|$; it is a positive $C^{\tau }
 $-function on $\overline{\uomega}$. Denote $J(y)^{\frac12}=\eta (y)$,
 then moreover,
$$
(Bu,v)_{L_2(\Omega )}=(\uB\,\uu,\eta
^2\bar{\uv})_{L_2(\uomega)}=(\eta \uB\eta ^{-1}\eta \uu,
\eta \uv)_{L_2(\uomega)};\tag6.6
$$
in particular one gets for $B=I$ that $ \|u\|_{L_2(\Omega)}=\|\eta \uu\|_{L_2(\uomega)}$.

Thus the mapping $\Phi :u\mapsto \eta \uu=\eta (\lambda ^{-1})^* u$ is an {\it isometry of}
$L_2(\Omega )$ {\it onto} $L_2(\uomega)$, i.e., a unitary operator
$\Phi =\eta (\lambda ^{-1})^*$, and by use of this mapping, $B$ is carried
over to the operator on $L_2(\uomega)$ $$
\widetilde B=\eta \uB\,\eta ^{-1}= \Phi B\Phi ^{-1}.\tag6.7$$
The analysis of the spectrum of $B$ can then be performed as an analysis
of the spectrum of $\widetilde B$. In particular, when $B$ is compact,
$s_j(B)=s_j(\widetilde B)$ for all
$j$.

We shall now show that the results of Section 5 extend to the general
case of \cite{AGW12}, where also the boundary is nonsmooth. 

\proclaim{Theorem 6.2} 
Theorems {\rm 5.2} and {\rm 5.4} (with corollaries and remarks) extend to the
situation where $\Omega $ is a bounded open set in ${\Bbb R}^n$ with
$B^{\frac32}_{p,2}$-boundary, satisfying {\rm (6.1)}.
\endproclaim

\demo{Proof} Consider the Dirichlet and Neumann-type problems for $A$
on $\Omega $, as defined in \cite{AGW12}. For each $l$,
apply the diffeomorphism $\lambda _l$ of Proposition 6.1 $2^\circ$ to carry the operators over to
$\uomega=\Omega _l$; as shown in \cite{AGW12} this gives operators
$A_l, C_l$ of the kind
treated in Section 5. For the resulting $\psi $do $Q_{l,+}$ and singular
Green terms $G_{\gamma ,l}=-K_{\gamma ,l}\gamma
_0Q_{l,+}$ and
$G_{C_l}=K_{\gamma ,l}L_l^{-1}(K_{\gamma ,l}')^*$ from the solution operators, the spectral asymptotics results of
Section 5 are valid.

Consider, for example, $B=G_C$ on $\Omega $, carried over to $\uB_l=G_{C_l }$ on $\Omega _l$ (the 
operators $Q_+$ and $G_\gamma $ are treated similarly). Theorem 5.4 and Remark 5.5 applied
to $\uB_l$ give that 
$$
\aligned
&s _j( \uB_l)j^{2/(n-1)}\to c( g^0_{C_l })^{2/(n-1)} \text{ for
}j\to\infty, \text{ where} \\
& c(g^0_{C_l })=\tfrac1{(n-1) (2\pi )^{(n-1)}}\int_{\Sigma
 _l}\int_{|\xi '|=1}(4(l_l^0)^2\operatorname{Re}\kappa
_l^+\operatorname{Re}\kappa _l^-)^{-(n-1)/4} \, d\omega dx'.
\endaligned \tag6.8
$$
To obtain similar results for $B$ we define as in (6.7):
$$
\widetilde B_l=\eta _l\uB_l\eta _l^{-1};\tag6.9
$$
here the $s$-numbers satisfy for all $j,l$,
$$
s_j(B)=s_j(\widetilde B_l).\tag6.10
$$
The results of Section 5 will not be applied directly to the operator $\eta
_l\uB_l\eta _l^{-1}$, since the multiplication by the $C^{\tau
_0}$-functions $\eta _l,\eta _l^{-1}$ does not fit easily into the
calculus. Instead we shall use their properties for $l\to\infty $.
Since $\lambda _l$ converges to the identity in $C^{1+\tau _0}$ for
$l\to\infty $ on a neighborhood $\Omega '$ of $\comega$, the Jacobian
$J_l$ converges to 1 in $C^{\tau _0}(\Omega' )$, and
so does its square root $\eta _l$ as well as the inverse function $\eta _l^{-1}$.

Now by (2.4), noting also that $\uB_l=\eta _l^{-1}\widetilde B_l\eta _l$,
$$
\|\eta _l\|_{L_\infty }s_j(\uB_l)\|\eta _l^{-1}\|_{L_\infty }\ge s_j(\widetilde B_l)\ge \frac1{\|\eta _l^{-1}\|_{L_\infty }}s_j(\uB_l)\frac1{\|\eta _l\|_{L_\infty }},
$$
 for all $j,l$ (norms in $L_\infty (\Omega ')$). Hence, in view of (6.10),
$$
\aligned
{\lim\sup}_{j\to\infty }s_j(B)j^{2/(n-1)}&\le \|\eta _l\|_{L_\infty }\|\eta _l^{-1}\|_{L_\infty }c(
g^0_{C_l })^{2/(n-1)} ,\\
{\lim\inf}_{j\to\infty }s_j(B)j^{2/(n-1)}&\ge \frac1{\|\eta _l\|_{L_\infty }\|\eta _l^{-1}\|_{L_\infty }}c(
g^0_{C_l })^{2/(n-1)} ,
\endaligned\tag6.11
$$
for all $l$. Since $\|\eta _l\|_{L_\infty }$ and $\|\eta _l^{-1}\|_{L_\infty }$ converge to 1
and $c(g^0_{C_l })^{2/(n-1)} \to c(
g^0_{C})^{2/(n-1)} $ for $l\to \infty $, 
we conclude that
$$
{\lim}_{j\to\infty }s _j( B)j^{2/(n-1)}= c( g^0_{C})^{2/(n-1)},\tag6.12
$$
as was to be shown.\qed
\enddemo

\head {Appendix. Proof of Theorem 3.7} \endhead
 
\head By Helmut Abels \endhead

For the following let $(\varphi_j(\xi))_{j\in\N_0}$ be a smooth partition of unity on $\Rn$ such that
$$
  \supp\{ \varphi_0\}\subseteq \ol{B_2(0)},\qquad  \supp\{ \varphi_0\}\subseteq \{\xi\in\Rn: 2^{j-1}\leq |\xi|\leq 2^{j+1}\}\quad \text{for}\ j\in \N
$$
satisfying $\varphi_j(-\xi)=\varphi_j(\xi)$ for all $\xi\in\Rn$, $j\in\N_0$ and
$$
  |\partial_\xi^\alpha \varphi_j(\xi)|\leq C_\alpha 2^{-j|\alpha|}\qquad \text{for all}\ j\in\N_0
$$
for arbitrary $\alpha\in\N_0^n$. This implies
$$ 
  |\partial_\xi^\alpha \varphi_j(\xi)|\leq C_\alpha
   \weight{\xi}^{-|\alpha|}\quad  \text{uniformly in}\ j\in\N_0\tag A.1 
$$
for all $\alpha \in \N_0^n$,
since $2^{j-1}\leq|\xi|\leq 2^{j+1}$ on $\supp \varphi_j$ if $j\geq 1$ and the
estimate for $j=0$ is trivial.

\proclaim{Lemma A.1}
  Let $Z$ be a Banach space, $N\in\N_0$, and let $g\colon \Rn\to Z$ be an $N$-times continuously differentiable function, $N\in \N_0$, with compact support.
  Then
  $$
    \|\F^{-1}[g](x)\|_Z\leq C_N|\supp g||x|^{-N}\sup_{|\beta|=
N}\|\partial_\xi^\beta g\|_{L_\infty(\Rn;Z)},\tag A.2
  $$
  uniformly in $x\neq 0$ and $g$.
\endproclaim

\demo{Proof}
  See e.g.\ Lemma 5.10 in Abels \cite{A12}. 
\qed\enddemo

\proclaim{Lemma A.2}
  Let $X_0,X_1$ be Banach spaces, let $p\in C^\tau S^m_{1,0}(\Rn\times\Rn,N;\Cal{L}(X_0,X_1))$, $m\in\R$, $N\in\N$, $\tau>0$, and let $p_j(x,\xi)=
  p(x,\xi)\varphi_j(\xi)$, $j\in\N_0$. Then for all
  $\alpha\in\N_0^n,M\in\{0,\ldots, N\}$ there is a constant $C_{\alpha,M}$ such that
  $$ 
    \| \partial_z^\alpha k_j(x,z)\|_{\Cal{L}(X_0,X_1)}
    \leq C_{\alpha,M} |z|^{-M} 2^{j(n+m+|\alpha|-M)}\tag A.3
  $$
  uniformly in $j\in\N_0$, $z\neq 0$, where
  $$
    k_j(x,z):= \F^{-1}_{\xi\mapsto z}[p_j(x,\xi)]\qquad \text{for all}\ x,z\in\Rn.
  $$
\endproclaim

\demo{Proof}
  The lemma follows from Lemma A.1 
applied to $g_x(\xi)=
  (i\xi)^\alpha p_j(x,\xi)$, where $x\in\Rn$ is considered as
  a fixed parameter, cf. Lemma 5.14 in Abels \cite{A12}. Since some of the arguments will be modified below, we include the details. Because of (A.1) 
and the Leibniz formula, it is easy to show that for all $\alpha\in\N_0^n$ with $|\alpha|\leq N$,
  $$
    \|\partial_\xi^\alpha 
    p_j(x,\xi)\|_{\Cal{L}(X_0,X_1)}\leq C_{\alpha}\weight{\xi}^{m-|\alpha|}
  $$
  uniformly in $x,\xi\in\Rn$ and $j\in\N_0$. Moreover, this implies that for all $\alpha,\gamma\in\N_0^n$ with $|\gamma|\leq N$,
  $$
    \bigl\|\partial_\xi^\gamma \bigl(
    (i\xi)^\alpha p_j(x,\xi)\bigr)\bigr\|_{\Cal{L}(X_0,X_1)}\leq C_{\alpha}\weight{\xi}^{m+|\alpha|-|\gamma|}
  $$
  uniformly in $x,\xi\in\Rn$ and $j\in\N_0$, by another use of
  the Leibniz formula. Therefore Lemma A.1 
yields
  $$\aligned
    {\|\partial_z^\alpha k_j(x,z)\|_{\Cal{L}(X_0,X_1)}}
    & \leq  C_N |\supp p_j(x,\cdot)| |z|^{-M} \sup_{|\gamma|=M} \bigl\|\partial_\xi^\gamma \bigl(
    (i\xi)^\alpha p_j(x,\xi)\bigr)\bigr\|_{\Cal{L}(X_0,X_1)}\\
& \leq C_{\alpha,\beta,M} |z|^{-M}2^{j(n+m+|\alpha|-M)},
  \endaligned $$
  since $|\supp p_j(x,.)|\leq C 2^{jn}$ independently of $x\in\Rn$.
\qed\enddemo

Note that for the last lemma no smoothness in $x$ is needed; in fact boundedness in $x$ would suffice.

Recall the notation $p(x,D_x )=\operatorname{OP}(p(x,\xi ))$, cf.\ (2.9).
\proclaim{Lemma A.3}
  Let $X_0,X_1$ be Banach spaces and let $p\in C^\tau S^m_{1,0}(\Rn\times\Rn,n+1;\Cal{L}(X_0,X_1))$, $m\in\R$, $\tau>0$, and let $p_j(x,D_x)=
  p(x,D_x)\varphi_j(D_x)$. Then for any $1\leq q< \infty$,
  $$ 
    \|p_j(x,D_x)\|_{\Cal{L}(L_q(\Rn;X_0),L_q(\Rn;X_1))} \leq C2^{jm}
\quad \text{for all}\ j\in\N_0 ,\tag A.4
  $$
  where $C$ does not depend on $j$.
\endproclaim

\demo{Proof}
The proof is variant of the proof of Lemma 6.20 in Abels \cite{A12}. We include the details for the convenience of the reader and since some of its arguments will be modified below.

 First of all
  $$
    p_j(x,D_x) f = \int_{\Rn} k_j(x,x-y) f(y) dy \quad \text{for all}\ f\in \SD(\Rn;X_0),
  $${}
  where $k_j$ satisfies (A.3). 
According to these estimates
  $$
    \int_{\Rn} \|k_j(x,z) \|_{\Cal{L}(X_0,X_1)} dz\leq C \bigl(\int_{|z|\leq 2^{-j}} 2^{j(n+m)} dz +
     \int_{|z|> 2^{-j}} |z|^{-n-1}2^{j(m-1)} dz \bigr).
  $$
  The first comes from (A.3) 
for $M=0$ and the second comes from the
  choice $M=n+1$. Hence we get by a simple calculation
  $$
    \int_{\Rn} \|k_j(x,z) \|_{\Cal{L}(X_0,X_1)} dz \leq C2^{jm},
  $$
  which proves (A.4), 
since
  $$
    \left\|\int_{\Rn} k_j(x,x-y) f(y) dy\right\|_{L_q(\Rn;X_1)}\leq \sup_{x\in \Rn}\int_{\Rn}
    \|k_j(x,z)\|_{\Cal{L}(X_0,X_1)}dz \|f\|_{L_q(\Rn;X_0)}
  $$
for all $f\in \SD(\Rn;X_0)$ because of Young's inequality for convolution integrals.
\qed\enddemo

\proclaim{Lemma A.4}
  Let $X_0,X_1$ be Banach spaces and let $p\in C^\tau S^m_{1,0}(\Rn\times\Rn,n+1;\Cal{L}(X_0,X_1))$, $m\in\R$, $\tau>0$, and let $p_{j,k}(x,\xi)=
  \varphi_k(D_x)p(x,\xi)\varphi_j(\xi)$ for alle $j,k\in\N_0$. Then for any $1\leq q< \infty$,
  $$ 
    \|p_{j,k}(x,D_x)\|_{\Cal{L}(L_q(\Rn;X_0),L_q(\Rn;X_1))} \leq
C2^{-\tau k}2^{jm} \quad \text{for all}\ j,k\in\N_0 ,\tag A.5
  $$
  where $C$ does not depend on $j,k$.
\endproclaim

\demo{Proof}
First of all we have that
$$\aligned
  &{\sup_{x\in\Rn} \|\partial_\xi^\alpha \varphi_k(D_x) p(x,\xi)\|_{\Cal{L}(X_0,X_1)}}\\
 &\quad \leq C 2^{-\tau k} \|\partial_\xi^\alpha p(\cdot ,\xi)\|_{C^\tau(\Rn;\Cal{L}(X_0,X_1))}\leq C'\weight{\xi}^{m-|\alpha|}2^{-\tau k},
\endaligned $$
since $C^\tau (\Rn;Z)\subseteq B^\tau_{\infty,\infty}(\Rn;Z)$ with continuous embedding, for a general Banach space $Z$. If $Z=\C$, then the statement follows e.g. from Theorem 6.2.5 in \cite{BL76}. For a general Banach space $Z$ the proof directly carries over.
 Hence by the Leibniz formula
  $$\aligned
  &{\|\partial_{\xi}^\beta\bigl(\varphi_k(D_x)p(x,\xi)\varphi_j(\xi)\bigr)\|_{\Cal{L}(X_0,X_1)}}\\
&\quad\leq C2^{-\tau k} \|\partial_{\xi}^\beta \bigl( p(\cdot ,\xi)\varphi_j(\xi)\bigr)\|_{C^\tau(\Rn;\Cal{L}(X_0,X_1))} \leq C'\weight{\xi}^{m-|\alpha|}2^{-\tau k}
  \endaligned $$
uniformly in $j,k\in\N_0$, $\xi\in \Rn$.
Now, if $k_{j,k}(x,z)= \F^{-1}_{\xi\mapsto z}[ p_{j,k}(x,\xi)]$, then
  $$
    p_{j,k}(x,D_x) f = \int_{\Rn} k_{j,k}(x,x-y) f(y) dy \quad \text{for all}\ f\in \SD(\Rn;X_0).
  $$
Hence as in the proof of Lemma A.2 
one shows that
  $$
    \| \partial_z^\alpha k_{j,k}(x,z)\|_{\Cal{L}(X_0,X_1)}
    \leq C_{\alpha,M} |z|^{-M} 2^{j(n+m+|\alpha|-M)}2^{-\tau k}
  $$
  uniformly in $j,k\in\N_0$, $z\neq 0$, for all $|\alpha|=M\in
  \{0,\ldots, n+1\}$. Therefore the same estimates as in the proof of
  Lemma A.3 
(with an additional factor $2^{-\tau k}$) prove (A.5).
\qed\enddemo

\proclaim{Lemma A.5}
  Let $X_0,X_1$ be Banach spaces, 
let $p\in C^\tau S^0_{1,0}(\Rn\times\Rn,n+1;\Cal{L}(X_0,X_1))$, and
let $p_j(x,\xi)=p(x,\xi)\varphi_j(\xi)$. Then for any $1\leq q < \infty$, there is a constant $C$ such that
  $$ 
    \|\varphi_i(D_x) p_j(x,D_x)\|_{\Cal{L}(L_q(\Rn;X_0),L_q(\Rn;X_1))} \leq
    \cases
     C2^{-i\tau} &\text{if}\ i\geq j+4,\\
     C&\text{if}\ j-3\leq i\leq j+3,\\
     C2^{-j\tau} &\text{if}\ j\geq i+4,
    \endcases \tag A.6
  $$
  uniformly in $i,j\in\N_0$.
\endproclaim

\demo{Proof}
First of all, if $j-3\leq i\leq j+3$, then the estimate follows from (A.4).
To prove the other cases, we use that
$$
  \varphi_i(D_x) p_j(x,D_x)= \sum_{k=0}^\infty \varphi_i(D_x) p_{j,k}(x,D_x),
$$
where $p_{j,k}(x,\xi)= \varphi_k(D_x) p_j(x,\xi)$. Moreover,
$$\aligned 
 & {\int_{\Rn}p_{j,k}(x,D_x) f(x) \ol{g(x)}\,d x}
  = \int_{{\Bbb R}^{2n}}  e^{ix\cdot \xi} \varphi_k(D_x) p_j(x,\xi)\hat{f}(\xi) \ol{g(x)}\,\d\xi d x\\
 &\quad = \int_{{\Bbb R}^{2n}} p_j(x,\xi)\hat{f}(\xi) \ol{\varphi_k(D_x) e^{-ix\cdot \xi} g(x)}\,\d\xi d x\\
&\quad = \int_{{\Bbb R}^{3n}} e^{-i\eta \cdot x}p_j(x,\xi)\hat{f}(\xi)\varphi_k(\eta )  \ol{\hat{g}(\xi+\eta )}\,\d\xi\d \eta  d x\\
&\quad = \int_{{\Bbb R}^{3n}} e^{i(\xi-\xi')\cdot x}p(x,\xi)\varphi_j(\xi)\hat{f}(\xi) \varphi_k(\xi'-\xi)  \ol{\hat{g}(\xi')}\,\d\xi\d \xi' d x
\endaligned$$
for all $f\in \SD(\Rn;X_1), g\in \SD(\Rn)$.
Here for every $x\in\Rn$ the support of
$$
  \int_{\Rn }e^{i(\xi-\xi')\cdot x}p(x,\xi)\varphi_j(\xi)\hat{f}(\xi) \varphi_k(\xi'-\xi)  \ol{\hat{g}(\xi')}\d\xi
$$
with respect to $\xi'$ is contained in 
$$\aligned
 &{\{2^{k-1}\leq |\xi|\leq 2^{k+1}\}+\{2^{j-1}\leq |\xi|\leq 2^{j+1}\}}\\
 &\quad \subseteq \cases
   \{2^{k-2}\leq |\xi|\leq 2^{k+2}\} &\text{if}\ k\geq j+3,\\
   \{2^{j-2}\leq |\xi|\leq 2^{j+2}\} &\text{if}\ j\geq k+3,\\
   \{|\xi|\leq 2^{j+2}\} &\text{if}\ j\geq k,\\
 \{|\xi|\leq 2^{k+2}\} &\text{if}\ k\geq j.\\
  \endcases
\endaligned
$$
Now let us consider first the case $i\geq j+4$. Then
$$\aligned
  &{\int_{\Rn}\varphi_i (D_x)p_{j,k}(x,D_x) f(x) \ol{g(x)}\,d x}\\
&\quad= \int_{{\Bbb R}^{3n}} e^{i(\xi-\xi')\cdot
x}p(x,\xi)\varphi_j(\xi)\hat{f}(\xi) \varphi_k(\xi'-\xi)
\ol{\varphi_i(\xi')\hat{g}(\xi')}\d\xi\d \xi' d x = 0
\endaligned
$$
for all $f\in \SD(\Rn;X_1), g\in \SD(\Rn)$, if $|k-i|\geq 4$,
 since 
$$
\alignat{1}
  \{2^{k-2}\leq |\xi|\leq 2^{k+2}\}\cap \{2^{i-1}\leq |\xi|\leq 2^{i+1}\} &=\emptyset\qquad \text{if}\ k\geq i+4,\\
  \{|\xi|\leq 2^{k+2}\}\cap \{2^{i-1}\leq |\xi|\leq 2^{i+1}\} &=\emptyset\qquad \text{if}\ j\leq k\leq i-4,\\
  \{|\xi|\leq 2^{j+2}\}\cap \{2^{i-1}\leq |\xi|\leq 2^{i+1}\} &=\emptyset\qquad \text{if}\ k< j.
\endalignat
$$
Therefore
$$\aligned
  \varphi_i(D_x) p_j(x,D_x)= \sum_{k=i-3}^{i+3} \varphi_i(D_x) p_{j,k}(x,D_x),
\endaligned $$
which implies
$$\aligned
 &\|\varphi_i(D_x) p_j(x,D_x)\|_{\Cal{L}(L_q(\Rn;X_0),L_q(\Rn;X_1))}\\
&\quad \leq  \sum_{k=i-3}^{i+3}\|\varphi_i(D_x) p_{j,k}(x,D_x)\|_{\Cal{L}(L_q(\Rn;X_0),L_q(\Rn;X_1))}
\leq  C 2^{-\tau i}
\endaligned $$
due to (A.5). 

Finally, we consider the case $j\geq i+4$. 
$$\aligned
  &{\int_{\Rn}\varphi_i (D_x)p_{j,k}(x,D_x) f(x) \ol{g(x)}\,d x}\\
&\quad = \int_{{\Bbb R}^{3n}} e^{i(\xi-\xi')\cdot x}p(x,\xi)\varphi_j(\xi)\hat{f}(\xi) \varphi_k(\xi'-\xi)  \ol{\varphi_i(\xi')\hat{g}(\xi')}\,\d\xi\d \xi' d x= 0
\endaligned $$
for all $f\in \SD(\Rn;X_1), g\in \SD(\Rn)$
if $|k-j|\geq 4$, since 
$$
\alignat{1}
  \{2^{k-2}\leq |\xi|\leq 2^{k+2}\}\cap \{2^{i-1}\leq |\xi|\leq 2^{i+1}\} &=\emptyset\qquad \text{if}\ k\geq j+4,\\
  \{2^{j-2}\leq|\xi|\leq 2^{j+2}\}\cap \{2^{i-1}\leq |\xi|\leq 2^{i+1}\} &=\emptyset\qquad \text{if}\  k\leq j-4.
\endalignat
$$
Therefore
$$\aligned
  \varphi_i(D_x) p_j(x,D_x)= \sum_{k=i-3}^{i+3} \varphi_i(D_x) p_{j,k}(x,D_x),
\endaligned $$
which implies
$$\aligned
  &{\|\varphi_i(D_x) p_j(x,D_x)\|_{\Cal{L}(L_q(\Rn;X_0),L_q(\Rn;X_1))}}\\
&\quad\leq  \sum_{k=j-3}^{j+3}\|\varphi_i(D_x) p_{j,k}(x,D_x)\|_{\Cal{L}(L_q(\Rn;X_0),L_q(\Rn;X_1))}
\leq  C 2^{-\tau j}
\endaligned $$
due to (A.5) 
again.
\qed\enddemo

\demo{Proof of Theorem~3.7}
First of all, by a composition from the right with $\weight{D_x}^{-m}$ we can always reduce to the case $m=0$, which we consider in the following. 

We use that
$
  f = \sum_{j=0}^\infty  f_j,
$
where $f_j= \varphi_j(D_x) f$, $f\in\SD(\Rn;H_0)$. Since $\supp \varphi_j \cap
\supp \varphi_k=\emptyset$ if $|j-k|> 1$, 
$$\aligned
  p(x,D_x) f &= \sum_{j=0}^\infty p_j(x,D_x) f = \sum_{k=0}^\infty p_j(x,D_x)
  (f_{j-1} + f_j + f_{j+1})\\
 &= \sum_{j=0}^\infty p_j(x,D_x) \tilde{f}_j,
\endaligned $$
where $\tilde{f}_j= f_{j-1} + f_j + f_{j+1}$ and we have set $f_{-1}=0$.
Moreover, we use that $\|f\|_{H^s(\Rn;H)}$ is equivalent to 
$$
  \|f\|_{B^s_{2,2}(\Rn;H)} :=\|(2^{sj}f_j)_{j\in\N_0}\|_{\ell_2(\N_0,L_2(\Rn;H))}
$$
if $H$ is a Hilbert space. This follows easily from Plancherel's theorem and the fact that $\weight{\xi}^s$ is equivalent to $2^{js}$ on $\supp \varphi_j$. 

If $f\in \SD(\Rn;H_0)$, then 
$$\aligned
&{ \|(2^{sj}\tilde{f}_j)_{j\in\N_0}\|_{\ell_2(\N_0,L_2(\Rn;H_0))}}
 = 
\|(2^{sj}(f_{j-1}+f_j+f_{j+1}))_{k\in\N_0}\|_{\ell_2(\N_0,L_2(\Rn;H_0))}
\\
  &\quad \leq C\bigl\| \bigl(2^{sj}\varphi_j(D_x)f\bigr)_{j\in\N_0}\bigr\|_{\ell_2(\N_0,L_2(\Rn;H_0))} = C\|f\|_{B^{s}_{2,2}(\Rn;H_0)}. 
\endaligned $$
Now using 
$$\aligned
  &{2^{is} \|\varphi_i(D_x) p_j(x,D_x)\tilde{f}_j\|_{L_2(\Rn;H_1)}}\\
&\quad \leq
  \cases
  C 2^{-(\tau -s)i-js} 2^{js}\|\tilde{f}_j\|_{L_2(\Rn;H_0)}&\text{if}\ i\geq j+4,\\  C 2^{js}\|\tilde{f}_j\|_{L_2(\Rn;H_0)}&\text{if}\ j-3\leq i\leq j+3,\\
  C 2^{(i-j)s-\tau j} 2^{js}\|\tilde{f}_j\|_{L_2(\Rn;H_0)}&\text{if}\ j\geq i+4,
  \endcases\\  
  &\quad\leq  C2^{-(\tau-|s|)|i-j|} 2^{js}\|\tilde{f}_j\|_{L_2(\Rn;H_0)}
\endaligned $$
Now let $a_j= 2^{-(\tau-|s|)|j|}$, $b_j= 2^{sj}\|\tilde{f}_j\|_{L_2(\Rn;H_0)}$ for $j\in\N_0$
and $a_j=b_j =0$ for $j\in\Z\setminus\N_0$. Then
$$\aligned
&{\|p(x,D_x) f\|_{B^s_{2,2}(\Rn;H_1)}}\\
 &\quad\leq C\bigl\| (a_j)_{j\in \Z}\ast (b_j)_{j\in
  \Z}\bigr\|_{\ell_2(\Z)}\leq C\|(a_j)_{j\in \Z}\|_{\ell_1(\Z)} \|(b_j)_{j\in
  \Z}\|_{\ell_2(\Z)}\\
&\quad = C  \|(2^{sj}\tilde{f}_j)_{j\in\N_0}\|_{\ell_2(\N_0,L_2(\Rn;H_0)))} \leq
  C \|f\|_{B^s_{2,2}(\Rn;H_0)},
\endaligned $$
where the sequence $(c_j)_{j\in\Z}= (a_j)_{j\in \Z}\ast (b_j)_{j\in \Z}$
  is the convolution of $(a_j)_{j\in \Z}$ and $(b_j)_{j\in \Z}$. (Recall that $B^s_{2,2}=H^s$.) Since $\SD(\Rn;H_0)$ is dense in $H^s(\Rn;H_0)$, the theorem follows. 
\qed\enddemo

\head Acknowledgement \endhead

The main author is grateful to Helmut Abels for many useful
conversations, and in particular for his helpfulness in working out the proof of
Theorem 3.7 in the Appendix. We also thank Heiko Gimperlein for
fruitful discussions.

\enddocument